\documentclass[11pt]{article}
\usepackage{amsfonts}
\usepackage{amsmath,amsbsy}
\usepackage{epsfig}
\usepackage{graphics}

\newtheorem{theorem}{Theorem}
\newtheorem{corollary}{Corollary}
\newtheorem{lemma}{Lemma}

\begin{document}

\centerline{\bf \Large GLOBAL EXISTENCE RESULTS FOR COMPLEX}

\bigskip

\centerline{\bf \Large HYPERBOLIC MODELS OF BACTERIAL CHEMOTAXIS}

\medskip
\bigskip

\centerline{\large Radek Erban$^*$ and Hyung Ju Hwang$^\dagger$}

\bigskip \medskip

{
\centerline{$^*$University of Oxford, Mathematical Institute} %
\centerline{24-29 St Giles', Oxford, OX1 3LB, United Kingdom} %
\centerline{e-mail: \textit{erban@maths.ox.ac.uk}} \medskip %
\centerline{$^\dagger$Trinity College Dublin, School of Mathematics} %
\centerline{Dublin 2, Ireland} 
\centerline{e-mail:
\textit{hjhwang@maths.tcd.ie}} }

\bigskip \medskip

\begin{abstract} \noindent
Bacteria are able to respond to environmental signals by changing their
rules of movement. When we take into account chemical signals in the
environment, this behaviour is often called chemotaxis. At the
individual-level, chemotaxis consists of several steps. First, the cell
detects the extracellular signal using receptors on its membrane. Then, the
cell processes the signal information through the intracellular signal
transduction network, and finally it responds by altering its motile
behaviour accordingly. At the population level, chemotaxis can lead to
aggregation of bacteria, travelling waves or pattern formation, and the
important task is to explain the population-level behaviour in terms of
individual-based models. It has been previously shown that the transport
equation framework \cite{Erban:2004:ICB,Erban:2004:STS} is suitable for
connecting different levels of modelling of bacterial chemotaxis. In this
paper, we couple the transport equation for bacteria with the
(parabolic/elliptic) equation for the extracellular signals. We prove global
existence of solutions for the general hyperbolic chemotaxis models of cells
which process the information about the extracellular signal through the
intracellular biochemical network and interact by altering the extracellular
signal as well. The conditions for global existence in terms of the
properties of the signal transduction model are given.
\end{abstract}

\section{Introduction}

\label{secintro}

The flagellated bacteria (e.g. \textit{Escherichia coli}, \textit{%
Salmone\-lla typhimurium}, \textit{Bacillus subtilis}) are single-celled
organisms. They are usually too small to be visible by the naked eye;
typically, they have the size of microns (see \cite%
{Salyers:2001:MDD,McKane:1996:MEA} for review). The behaviour of a bacterium
is primarily influenced by concentrations of various chemicals inside the
cell. Since bacteria are small, we can assume that the concentrations of the
chemicals inside the cytoplasm are uniform. Therefore, we can suppose that
the cells are points. Moreover, to create a mathematical description of a
bacterium, we introduce the vector of internal state variables \cite%
{Erban:2004:ICB,Erban:2004:STS,Erban:2005:ICB} 
\begin{equation}
\mathbf{y}=(y_{1},y_{2},\dots ,y_{m})^{T}\in \mathbb{R}^{m},  \label{rom200}
\end{equation}%
where $y_{i},$ $i=1,\dots ,m,$ are concentrations of various chemicals
(proteins, receptor states etc.) inside the cell involved in the processes
of interest. The individual behaviour of a cell primarily depends on the
vector $\mathbf{y}$ which is a function of time. Consequently, the state of
bacterium is uniquely determined by the vector $(t,\mathbf{x},\mathbf{v},%
\mathbf{y})$ where $\mathbf{x}\in \mathbb{R}^{N}$ is the position of a cell, 
$\mathbf{v}\in \mathbb{R}^{N}$ is its velocity, $\mathbf{y}\in \mathbb{R}^{m}
$ is its internal state, $t$ is time and $N=1,2,$ or $3,$ is the dimension
of the physical space.

During its life, a cell must communicate with its environment in order to
find nutrients, to avoid repellents, to find mates etc. For this purpose,
there are receptors in the cellular membrane which can detect various
chemicals in the environment. We describe the chemicals outside the
bacterium by the signaling vector (which depends on the position of the cell 
$\mathbf{x}$ and time $t$) 
\begin{equation}
\mathbf{S}(\mathbf{x},t)=(S_{1},S_{2},\dots ,S_{M})^T \in \mathbb{R}^{M}.
\label{rom201}
\end{equation}%
Then the evolution of the internal state vector $\mathbf{y}$ depends also on
the signaling vector $\mathbf{S}.$ Since we describe chemical processes, we
can assume that $\mathbf{y}$ evolves according to system of ordinary
differential equations 
\begin{equation}
\frac{\mbox{d}\mathbf{y}}{\mbox{d}t}=\mathbf{F}(\mathbf{S}(\mathbf{x}),%
\mathbf{y}).  \label{rom14}
\end{equation}%
This system formally captures all biochemistry inside the cell and
therefore, the concrete form of the vector function $\mathbf{F}:\mathbb{R}%
^{M}\times \mathbb{R}^{m}\rightarrow \mathbb{R}^{m}$ can be very complicated
depending on the number of details which are included in the model.

Bacterial movement and the signal transduction network (\ref{rom14}) will be
discussed in more details in Section \ref{secbiodetails}. From the
mathematical point of view, the movement of the flagellated bacteria can be
viewed as a biased random walk. The properties of this random walk depend on
the internal state $\mathbf{y}$ and bacterial velocity $\mathbf{v}$. The
classical description of the bacterial movement is the so called velocity
jump process \cite{Othmer:1988:MDB,Erban:2004:ICB,Erban:2004:STS}. It means
that the bacterium runs with some velocity and at random instants of time it
changes its velocity according to the Poisson process with the intensity $%
\lambda(\mathbf{y})$.

Let $f(\mathbf{x},\mathbf{v},\mathbf{y},t)$ be the density function of
bacteria in a $(2N+m)-$dimensional phase space with coordinates $(\mathbf{x},%
\mathbf{v},\mathbf{y})$ where $\mathbf{x}\in \mathbb{R}^{N}$ is the position
of a cell, $\mathbf{v}\in V\subset \mathbb{R}^{N}$ is its velocity and $%
\mathbf{y}\in \mathbb{R}^{m}$ is its internal state, which evolves according
to (\ref{rom14}). Thus $f(\mathbf{x},\mathbf{v},\mathbf{y},t)\mbox{d}\mathbf{%
x}\mbox{d}\mathbf{v}\mbox{d}\mathbf{y}$ is the number of cells with position
between $\mathbf{x}$ and $\mathbf{x}+\mbox{d}\mathbf{x}$, velocity between $%
\mathbf{v}$ and $\mathbf{v}+\mbox{d}\mathbf{v}$, and internal state between $%
\mathbf{y}$ and $\mathbf{y}+\mbox{d}\mathbf{y}.$ Then evolution of $f$ is
governed by the following transport equation \cite%
{Erban:2004:ICB,Erban:2004:STS} 
\begin{equation}
\frac{\partial f}{\partial t}+\nabla _{\mathbf{x}}\cdot \mathbf{v}f+\nabla _{%
\mathbf{y}}\cdot \mathbf{F}(\mathbf{S}(\mathbf{x}),\mathbf{y})f=-\lambda (%
\mathbf{y})f+\int_{V}\lambda (\mathbf{y})K(\mathbf{v},\mathbf{v}^{\prime },%
\mathbf{y})f(\mathbf{x},\mathbf{v}^{\prime },\mathbf{y},t)\mbox{d}\mathbf{v}%
^{\prime }  \label{vjpint}
\end{equation}%
where the kernel $K(\mathbf{v},\mathbf{v}^{\prime },\mathbf{y})$ gives the
probability of a change in velocity from $\mathbf{v}^{\prime }$ to $\mathbf{v%
}$, given that a reorientation occurs. We assume that the random velocity
changes are the result of a Poisson process of intensity $\lambda (\mathbf{y}%
)$ \cite{Block:1983:AKB}. The kernel $K$ is non-negative and satisfies the
normalization condition 
\begin{equation}
\int_{V}K(\mathbf{v},\mathbf{v}^{\prime },\mathbf{y})\mbox{d}\mathbf{v}=1
\label{normalcond}
\end{equation}%
where 
\begin{equation}
V\text{ \ is a symmetric compact set in }\mathbb{R}^{N}.  \label{setvel}
\end{equation}%
Realistic examples of the kernel $K(\mathbf{v},\mathbf{v}^{\prime },\mathbf{y%
})$, set $V$, signal transduction network $\mathbf{F}$ and the turning
frequency $\lambda (\mathbf{y})$ are given in Section \ref{secbiodetails}.
They all satisfy the above basic assumptions.

To write equation (\ref{vjpint}) in more compact form, we introduce the
kernel $T$ defined as a product of the turning frequency $\lambda $ and the
kernel $K$, i.e. 
\begin{equation}
T(\mathbf{v},\mathbf{v}^{\prime },\mathbf{y})=\lambda \left( \mathbf{y}%
\right) K(\mathbf{v},\mathbf{v}^{\prime },\mathbf{y}).  \label{kernelT}
\end{equation}%
Moreover, our goal is to couple equation (\ref{vjpint}) with the realistic
system of partial differential equations for the extracellular signal vector 
$\mathbf{S}.$ We assume that the external signal diffuses. It can be also
produced or degraded by bacteria, degraded on its own or the components of $%
\mathbf{S}$ can react with each other in the extracellular space. Hence, the
general hyperbolic system of interest can be written in the following form: 
\begin{eqnarray}
\frac{\partial f}{\partial t}+\nabla _{\mathbf{x}}\cdot \mathbf{v}f+\nabla _{%
\mathbf{y}}\cdot \mathbf{F}(\mathbf{S}(\mathbf{x}),\mathbf{y})f &=&\int_{V}T(%
\mathbf{v},\mathbf{v}^{\prime },\mathbf{y})\Big[f(\mathbf{v}^{\prime })-f(%
\mathbf{v})\Big]\mbox{d}\mathbf{v}^{\prime }  \label{hyperp} \\
\frac{\partial \mathbf{S}}{\partial t} &=&D\triangle \mathbf{S}+\mathbf{R}(%
\mathbf{S},n)  \label{parabS}
\end{eqnarray}%
where $n\equiv n(\mathbf{x},t)$ is the macroscopic density of individuals at
point $\mathbf{x}\in \mathbb{R}^{N}$ and time $t$ given as 
\begin{equation}
n(\mathbf{x},t)=\int_{\mathbb{R}^{m}}\int_{V}f(\mathbf{x},\mathbf{v},\mathbf{%
y},t)\mbox{d}\mathbf{v}\mbox{d}\mathbf{y},  \label{equationforn}
\end{equation}%
D is a diagonal $M\times M$ matrix which diagonal elements are diffusion
constants of different chemicals in the extracellular signal vector $\mathbf{%
S}$ and the term $\mathbf{R}:\mathbb{R}^{M}\times \mathbb{R}\rightarrow 
\mathbb{R}^{M}$ describes the creation, reaction and degradation of the
signals.

The goal of this paper is to prove global existence results for the system (%
\ref{hyperp}) -- (\ref{parabS}). We will focus on one-dimensional case in
what follows. In Section \ref{secmotiv}, we will start with a simple model
of the signal transduction (\ref{simmodel}) which was used previously \cite%
{Erban:2004:ICB,Erban:2004:STS}. The simple model (\ref{simmodel}) has the
essential properties of the realistic models of the signal transduction, but
it is more tractable from the mathematical point of view than more complex
models of bacterial chemotaxis. We prove the global existence of solutions
of the one-dimensional version of (\ref{hyperp}) -- (\ref{parabS}) with the
simplified model of the signal transduction.

In order to study the general case, we first review the relevant biology in
Section \ref{secbiodetails}. This will help us to specify the realistic
conditions on signal transduction model $(\ref{rom14})$, turning frequency $%
\lambda(\mathbf{y})$, turning kernel $K(\mathbf{v},\mathbf{v}^{\prime },%
\mathbf{y})$, set $V$, diffusion matrix $D$ and the reaction term $\mathbf{R}%
(\mathbf{S},n)$ in equations (\ref{hyperp}) --(\ref{parabS}). In Section \ref%
{secglobalexistence}, we study the global existence for the system (\ref%
{hyperp}) -- (\ref{parabS}) for general models of bacterial signal
transduction which are introduced in Section \ref{secbiodetails}. We also
consider that equation (\ref{parabS}) is at quasi-equilibrium, i.e. we
consider the elliptic equation for the signal in Section \ref%
{secglobalexistence}.

Hence, this paper consists of two main mathematical results. First, we prove
the global existence of solutions to the problem (\ref{hyperp}) -- (\ref%
{parabS}) for the simplified model of signal transduction (\ref{simmodel})
and for the system of parabolic equations (\ref{parabS}) for the
extracellular signal (see Section \ref{secmotiv}, Theorem \ref{theorem1}).
Then, we prove the global existence of solutions for the general model of
signal transduction (\ref{rom14}) coupled with the system of elliptic
equations for the extracellular signal (see Section \ref{secglobalexistence}%
, Theorem \ref{theorem2}). The necessary growth assumptions on turning
frequency $T(\mathbf{v},\mathbf{v}^{\prime },\mathbf{y})$ are given in terms
of the signal derivative along the cell trajectory. It means that the growth
estimates on $T$ include the temporal derivative as well as the spatial
derivative of the extracellular signal. Finally, we provide discussion and
comparison with relevant results from the literature in Section \ref%
{secdiscussion}.

\section{Global existence for a simplified model of signal transduction}

\label{secmotiv}

A simplified model of excitation-adaptation dynamics was studied in \cite%
{Erban:2004:ICB,Erban:2004:STS,Othmer:1998:OCS,Dolak:2005:KMC} 
where $\mathbf{y}%
=(y_{1},y_{2})^{T}\in \mathbb{R}^{2}$ and the right hand side of equation (%
\ref{rom14}) was given as 
\begin{equation}
\mathbf{F}\equiv \left( 
\begin{array}{c}
F_{1} \\ 
F_{2}%
\end{array}%
\right) =\left( 
\begin{array}{c}
\displaystyle\frac{g(\mathbf{S}(\mathbf{x},t))-(y_{1}+y_{2})}{t_{e}}%
\raisebox{-5.6mm}{\rule{0pt}{12.2mm}} \\ 
\displaystyle\frac{g(\mathbf{S}(\mathbf{x},t))-y_{2}}{t_{a}}%
\end{array}%
\right)   \label{simmodel}
\end{equation}%
where $t_{e}$ and $t_{a}$ are positive constants and $g:\mathbb{R}%
^{M}\rightarrow \lbrack 0,\infty )$. We will see in Section \ref%
{secbiodetails} that the simplified model (\ref{simmodel}) has the essential
properties of realistic signal transduction models. Hence, the model (\ref%
{simmodel}) is a natural starting point of this paper. For simplicity, we
work in a  one-dimensional physical space, i.e. $N=1,$ and the goal of this
section is to prove Theorem \ref{theorem1} about the system (\ref{hyperp})
-- (\ref{parabS}). In what follows, we denote $L^{p}(\Omega )$, $1\leq p\leq
\infty $, $\Omega \subset \mathbb{R}^{d}$, the Banach space of measurable
functions with the finite norms 
\begin{equation*}
\Vert h\Vert \hbox{\raise -0.5mm \hbox{$_{L^p}$}}=\left( \int_{\Omega }|h(%
\mathbf{x})|^{p}\mbox{d}\mathbf{x}\right) ^{1/p},\;\;\mbox{for}\;1\leq
p<\infty ,\qquad \mbox{and}\qquad \Vert h\Vert 
\hbox{\raise -0.5mm
\hbox{$_{L^\infty}$}}=\mathop{\mbox{ess}\,\mbox{sup}}_{\Omega }|h(\mathbf{x}%
)|.
\end{equation*}%
We denote $W^{k,p}(\Omega )$, $1\leq p\leq \infty $, $\Omega \subset \mathbb{%
R}^{d}$, the usual Sobolev space 
\begin{equation*}
W^{k,p}(\Omega )=\left\{ h\in L^{p}(\Omega )\,|\,\forall \boldsymbol{\alpha}%
\in \mathbb{N}_{0}^{d},|\boldsymbol{\alpha}|\leq k\Rightarrow \frac{\partial
^{|\boldsymbol{\alpha}|}h}{\partial x_{1}^{\alpha _{1}}x_{2}^{\alpha
_{2}}\dots x_{d}^{\alpha _{d}}}\in L^{p}(\Omega )\right\} 
\end{equation*}%
where $\boldsymbol{\alpha}=(\alpha _{1},\alpha _{2},\dots ,\alpha _{d})\in 
\mathbb{N}_{0}^{d}$ is a vector of nonnegative integers and $|%
\boldsymbol{\alpha}|=\alpha _{1}+\alpha _{2}+\dots +\alpha _{d}.$ The norm
in $W^{k,p}(\Omega )$ is defined as 
\begin{equation*}
\Vert h\Vert \hbox{\raise -0.5mm \hbox{$_{W^{k,p}}$}}=\sum_{%
\boldsymbol{\alpha}\in \mathbb{N}_{0}^{d},\,|\boldsymbol{\alpha}|\leq
k}\left\Vert \frac{\partial^{|\boldsymbol{\alpha}|}h}{\partial
x_{1}^{\alpha _{1}}x_{2}^{\alpha _{2}}\dots x_{d}^{\alpha_{d}}}\right\Vert
_{L^{p}}
\end{equation*}%
To simplify mathematical formulas, we will make use of the following
notation. For any function $h:\mathbb{R}^{d}\rightarrow \mathbb{R}$, $\nabla
h$ denotes the gradient of $h$ with respect to all variables and $\nabla
_{\!x_{1}x_{2}}h$ is the 2-dimensional gradient vector with respect to the
variables $x_{1}$ and $x_{2}$ only, i.e. 
\begin{equation}
\nabla h=\left( \frac{\partial h}{\partial x_{1}},\frac{\partial h}{\partial
x_{2}},\dots ,\frac{\partial h}{\partial x_{d}}\right) ,\qquad \mbox{and}%
\qquad \nabla _{\!x_{1}x_{2}}h=\left( \frac{\partial h}{\partial x_{1}},%
\frac{\partial h}{\partial x_{2}}\right) .  \label{gradnotation}
\end{equation}%
We already made use of this notation in equation (\ref{hyperp}) where the
gradients of the function $f$ were taken only with respect to the selected
parts of the state vector. In this section, we study the movement of cells
in one dimension, i.e. $N=1$. Moreover, we assume that the external signal
diffuses and it is produced by bacteria and degraded on its own. Hence, the
system of equations (\ref{hyperp}) -- (\ref{parabS}) reads as follows 
\begin{equation}
\frac{\partial f}{\partial t}+\nabla _{\!x}\cdot vf+\nabla _{\!\mathbf{y}%
}\cdot \mathbf{F}(\mathbf{S}(\mathbf{x}),\mathbf{y})f=\int_{V}T(v,v^{\prime
},\mathbf{y})\Big[f(v^{\prime })-f(v)\Big]\mbox{d}v^{\prime }
\label{hyperpsim}
\end{equation}%
\begin{equation}
\frac{\partial S_{i}}{\partial t}=d_{i}\frac{\partial ^{2}S_{i}}{\partial
x^{2}}+k_{i}n-k_{i}^{0}S_{i},\qquad \;i=1,\dots ,M,  \label{parabSsim}
\end{equation}%
where $d_{i}$, $k_{i}$ and $k_{i}^{0}$ are positive constants and $n\equiv
n(x,t)$ is the macroscopic density of individuals at point $x\in \mathbb{R}$
and time $t$ given by (\ref{equationforn}). Position $x$ and velocity $v$
are scalars for $N=1$, so we do not use bold letters for position and
velocity in equation (\ref{hyperpsim}). Otherwise, equation (\ref{hyperpsim}%
) is the same as equation (\ref{hyperp}). Following notation (\ref%
{gradnotation}), symbol $\nabla _{\!x}f$ denotes the partial derivative of
distribution function $f$ with respect to $x.$ Let us note that (depending
on the form of function $g$ in (\ref{simmodel})) some extracellular signals
might be attractants and some extracellular signals might be repellents. If
we have sufficient growth estimates on function $g$ and kernel $T$, we can
guarantee the global existence of solutions of system (\ref{hyperpsim}) -- (%
\ref{parabSsim}) as it is shown in the following theorem.

\begin{theorem}
Consider that the function $\mathbf{F}$ is given by $(\ref{simmodel})$.
Assume that there exist non-decreasing positive continuous functions $\Phi
,\Psi \in C\left( \mathbb{R}\right) $ satisfying 
\begin{equation}
\left\vert g(\mathbf{z})\right\vert +\left\vert \nabla g(\mathbf{z}%
)\right\vert \leq \Phi \left( \left\vert \mathbf{z}\right\vert \right) \quad 
\text{and}\quad \left\vert T(v,v^{\prime },\mathbf{y})\right\vert
+\left\vert \nabla T(v,v^{\prime },\mathbf{y})\right\vert \leq \Psi \left(
\left\vert \mathbf{y}\right\vert \right) .  \label{growthassumptions}
\end{equation}%
Assume that $f_{0}\in W^{1,1}({}\mathbb{R}\times V\times \mathbb{R}^{2})\cap
W^{1,\infty }({}\mathbb{R}\times V\times \mathbb{R}^{2})$ with compact
support and $\mathbf{S}_{0}\in \lbrack W^{1,\infty }({}\mathbb{R})]^{M}$
with compact support. Then there exist global solutions of the system $(\ref%
{hyperpsim})$ -- $(\ref{parabSsim})$ satisfying 
\begin{equation}
f(\cdot ,\cdot ,\cdot ,t)\in W^{1,1}({}\mathbb{R}\times V\times \mathbb{R}%
^{2})\cap W^{1,\infty }({}\mathbb{R}\times V\times \mathbb{R}^{2}),
\label{spacef}
\end{equation}%
\begin{equation}
\mathbf{S}(\cdot ,t)\in \left[ W^{1,\infty }(\mathbb{R})\right] ^{M}
\label{spaceS}
\end{equation}%
and initial conditions $f(\cdot ,\cdot ,\cdot ,0)=f_{0}(\cdot ,\cdot ,\cdot )
$ and $\mathbf{S}(\cdot ,0)=\mathbf{S}_{0}(\cdot )$. \label{theorem1}
\end{theorem}

\bigskip

\noindent First, the characteristics of the hyperbolic equation (\ref%
{hyperpsim}) are given for $N=1$ as 
\begin{equation}
\frac{\mbox{d}X}{\mbox{d}s}=V,\quad \frac{\mbox{d}V}{\mbox{d}s}=0,\quad 
\frac{\mbox{d}\mathbf{Y}}{\mbox{d}s}=\mathbf{F}(\mathbf{S}\left( X\left(
s\right) ,s\right) ,\mathbf{Y}(s)).  \label{chareqn}
\end{equation}%
Then along back-time characteristics starting at $\left( x,v,\mathbf{y}%
,t\right) $, we have for $0\leq s\leq t,$%
\begin{eqnarray}
X\left( s;x,v\mathbf{,y},t\right) &=&x-v\left( t-s\right) ,\ \   \label{X} \\
\mathbf{Y}\left( s;x,v\mathbf{,y},t\right) &=&\mathbf{y}-\int_{s}^{t}\mathbf{%
F}\left( \mathbf{S}\left( X(\tau ),\tau \right) ,\mathbf{Y}(\tau )\right) %
\mbox{d}\tau .  \label{Y}
\end{eqnarray}%
Next, we will prove several auxiliary lemmas.

\begin{lemma}
Derivation of the characteristics $(\ref{X})$ and $(\ref{Y})$ with respect
to the initial conditions gives, for $0\leq s\leq t,$ 
\begin{equation}
\frac{\partial X}{\partial x}=1,\quad \frac{\partial \mathbf{Y}}{\partial 
\mathbf{y}}=\exp \left[ \frac{\partial \mathbf{F}}{\partial \mathbf{y}}(s-t)%
\right] \qquad \mbox{where}\quad \frac{\partial \mathbf{F}}{\partial \mathbf{%
y}}=\left( 
\begin{array}{cc}
-\displaystyle\frac{1}{t_{e}}\raisebox{-5.6mm}{\rule{0pt}{12.2mm}} & -%
\displaystyle\frac{1}{t_{e}} \\ 
0 & -\displaystyle\frac{1}{t_{a}}%
\end{array}%
\right) .  \label{derchar}
\end{equation}%
Moreover, 
\begin{equation}
\det \frac{\partial \mathbf{Y}}{\partial \mathbf{y}}=\exp \left[ \left( 
\frac{1}{t_{e}}+\frac{1}{t_{a}}\right) \left( t-s\right) \right] \geq 1.
\label{determyy}
\end{equation}%
\label{lemma3}
\end{lemma}

\noindent \textbf{Proof.} We differentiate (\ref{Y}) with respect to $%
\mathbf{y}$ to get%
\begin{equation}
\frac{\partial \mathbf{Y}}{\partial \mathbf{y}}=\mathbf{I}_{2}+\int_{t}^{s}%
\frac{\partial \mathbf{F}}{\partial \mathbf{y}}\frac{\partial \mathbf{Y}}{%
\partial \mathbf{y}}\left( \tau \right) d\tau .  \label{derY}
\end{equation}%
where $\mathbf{I}_{2}$ is the $2\times 2$ identity matrix. Let 
\begin{equation*}
\mathbf{G}\left( s\right) =\int_{t}^{s}\frac{\partial \mathbf{F}}{\partial 
\mathbf{y}}\frac{\partial \mathbf{Y}}{\partial \mathbf{y}}\left( \tau
\right) d\tau ,
\end{equation*}%
then we have%
\begin{equation*}
\mathbf{G}^{\prime }(s)=\frac{\partial \mathbf{F}}{\partial \mathbf{y}}\frac{%
\partial \mathbf{Y}}{\partial \mathbf{y}}(s).
\end{equation*}%
Using (\ref{derY}), we obtain 
\begin{equation*}
\mathbf{G}^{\prime }\left( s\right) -\frac{\partial \mathbf{F}}{\partial 
\mathbf{y}}\mathbf{G}\left( s\right) =\frac{\partial \mathbf{F}}{\partial 
\mathbf{y}}.
\end{equation*}%
Integrating the last equation, we have 
\begin{equation*}
\mathbf{G}\left( s\right) =\exp \left[ \frac{\partial \mathbf{F}}{\partial 
\mathbf{y}}(s-t)\right] -\mathbf{I}_{2},
\end{equation*}%
which deduce (\ref{derchar}). Computing the determinant of (\ref{derchar}),
we derive (\ref{determyy}).

\rightline{Q.E.D.}

\bigskip

\begin{lemma}
Let us assume $(\ref{simmodel})$ and $(\ref{growthassumptions})$. Then the
solution of $(\ref{chareqn})$ satisfies 
\begin{equation}
|\mathbf{Y}(\tau)|\leq C\left\{ 1+\Phi \left( \sup_{0\leq s\leq \tau
}\left\vert \mathbf{S}\left( X\left( s\right) \!,s\right) \right\vert
\right) \right\}   \label{grad Y-S}
\end{equation}%
where $C$ depends on the $y$-support of $f_{0}$ and $\mathbf{S}_{0},$ $%
t_{a}, $ and $t_{e}.$ \label{lemma4}
\end{lemma}

\noindent \textbf{Proof.} Using the assumption (\ref{growthassumptions}) and
applying the Gronwall inequality to the ordinary differential equation (\ref%
{chareqn}) yields%
\begin{eqnarray*}
Y_{2}\left( \tau \right)  &=&Y_{2}\left( 0\right) \exp \left( -\tau
/t_{a}\right) +\frac{1}{t_{a}}\int_{0}^{\tau }g\left( \mathbf{S}(X\left(
s\right) \!,s)\right) \exp [(s-\tau )/t_{a}]ds \\
&\leq &\left\vert Y_{2}\left( 0\right) \right\vert +\frac{1}{t_{a}^{2}}\Phi
\left( \sup_{0\leq s\leq \tau }\left\vert \mathbf{S}\left( X\left( s\right)
\!,s\right) \right\vert \right) .
\end{eqnarray*}%
In a similar way, we get%
\begin{eqnarray*}
Y_{1}\left( \tau \right)  &=&Y_{1}\left( 0\right) \exp \left( -\tau
/t_{e}\right) +\frac{1}{t_{e}}\int_{0}^{\tau }\{g\left( \mathbf{S}(X\left(
s\right) \!,s)\right) -Y_{2}\left( s\right) \}\exp [(s-\tau )/t_{e}]ds \\
&\leq &\left\vert Y_{1}\left( 0\right) \right\vert +\frac{1}{t_{e}^{2}}%
\sup_{0\leq s\leq \tau }\left\vert Y_{2}\left( s\right) \right\vert +\frac{1%
}{t_{e}^{2}}\Phi \left( \sup_{0\leq s\leq \tau }\left\vert \mathbf{S}\left(
X\left( s\right) \!,s\right) \right\vert \right) .
\end{eqnarray*}%
Thus we deduce (\ref{grad Y-S}).

\rightline{Q.E.D.}

\begin{lemma}
\label{S-gradS} If $n\in L^{\infty }([0,\infty ):L^{1}\left( \mathbb{R}%
\right) \cap L^{2}\left( \mathbb{R}\right) ),$ then the solution $\mathbf{S}$
of the system of equations $(\ref{parabSsim})$ satisfies 
\begin{equation*}
\left\Vert \mathbf{S}\left( t\right) \right\Vert _{L^{\infty }}\leq
C\sup_{0\leq \tau \leq t}\left\Vert n\left( t\right) \right\Vert
_{L^{1}}=C\left\Vert n\left( 0\right) \right\Vert _{L^{1}},
\end{equation*}%
\begin{equation*}
\left\Vert \frac{\partial \mathbf{S}}{\partial x}\left( t\right) \right\Vert
_{L^{\infty }}\leq C\left[ 1+\left\Vert n\left( 0\right) \right\Vert
_{L^{1}}\left( 1+\left( \ln t\right) _{+}+\left\vert \ln \left( \sup_{0\leq
\tau \leq t}\left\Vert n\left( \tau \right) \right\Vert _{L^{2}}\right)
\right\vert \right) \right] 
\end{equation*}%
where $\left( \cdot \right) _{+}$ means the positive part and the constant $C
$ depends only on $k_{i}$, $k_{i}^{0}$ and $d_{i}.$ \label{lemma2}
\end{lemma}

\noindent\textbf{Proof.} See \cite[Lemma 4]{Hwang:2005:GEC}.

\rightline{Q.E.D.}

\bigskip

\noindent \textbf{Proof of Theorem \ref{theorem1}.} Integrating (\ref%
{hyperpsim}) along the characteristic (\ref{X}) -- (\ref{Y}) from $0$ to $t$
and using (\ref{growthassumptions}), we get 
\begin{equation*}
f\left( x,v,\mathbf{y},t\right) \leq f_{0}\left( X\left( 0\right) ,v,\mathbf{%
Y}\left( 0\right) \right) +
\end{equation*}%
\begin{equation*}
+\;C\int_{0}^{t}\Psi \left( \left\vert \mathbf{Y}(\tau )\right\vert \right)
\times \left[ f\left( X(\tau ),v,\mathbf{Y}(\tau ),\tau \right)
+\int_{V}f\left( X(\tau ),v^{\prime },\mathbf{Y}(\tau ),\tau \right)
dv^{\prime }\right] \mbox{d}\tau +
\end{equation*}%
\begin{equation*}
+\int_{0}^{t}\left\vert \nabla _{\mathbf{y}}\cdot \mathbf{F}\left( \mathbf{S(%
}X\left( \tau \right) ),\mathbf{Y}\left( \tau \right) \right) \right\vert
f\left( X\left( \tau \right) ,v,\mathbf{Y}\left( \tau \right) ,\tau \right) %
\mbox{d}\tau .
\end{equation*}%
Since $\nabla _{y}\cdot \mathbf{F=-}\frac{1}{t_{e}}-\frac{1}{t_{a}}$, we get
(using Lemma \ref{lemma4}) 
\begin{equation}
f\left( x,v,\mathbf{y},t\right) \;\leq \;f_{0}\left( X\left( 0\right) ,v,%
\mathbf{Y}\left( 0\right) \right) +C\int_{0}^{t}f\left( X\left( \tau \right)
,v,\mathbf{Y}\left( \tau \right) ,\tau \right) \mbox{d}\tau +
\label{ineqpom}
\end{equation}%
\begin{equation*}
+C\int_{0}^{t}\Psi \left( C\left[ 1+\Phi \left( \sup_{0\leq s\leq \tau
}\left\vert \mathbf{S}\left( X\left( s\right) ,s\right) \right\vert \right) %
\right] \right) \times 
\end{equation*}%
\begin{equation*}
\times \left[ |V|f\left( X(\tau ),v,\mathbf{Y}(\tau ),\tau \right)
+\int_{V}f\left( X(\tau ),v^{\prime },\mathbf{Y}(\tau ),\tau \right)
dv^{\prime }\right] \mbox{d}\tau 
\end{equation*}%
where $C$ is a constant depending only on support of $f_{0}$, $S_{0}$, $t_{e}
$ and $t_{a}.$ Using Lemma \ref{lemma3}, we have 
\begin{equation}
\left( \det \frac{\partial \mathbf{Y}}{\partial \mathbf{y}}\right) ^{-1}\leq
1,\qquad \left( \det \frac{\partial X}{\partial x}\right) ^{-1}=1
\label{det}
\end{equation}%
and so%
\begin{equation*}
\int_{\mathbb{R}\times V\times \mathbb{R}^{M}}\int_{V}f^{p}\left( X\left(
\tau \right) ,v^{\prime },\mathbf{Y}\left( \tau \right) ,\tau \right)
dv^{\prime }\mbox{d}x\mbox{d}v\mbox{d}\mathbf{y}=
\end{equation*}%
\begin{equation*}
=\left\vert V\right\vert \int f^{p}\left( X\left( \tau \right) ,v^{\prime },%
\mathbf{Y}\left( \tau \right) ,\tau \right) \left( \det \frac{\partial 
\mathbf{Y}}{\partial \mathbf{y}}\right) ^{-1}\left( \det \frac{\partial X}{%
\partial x}\right) ^{-1}dv^{\prime }\mbox{d}X\mbox{d}\mathbf{Y}\leq 
\end{equation*}%
\begin{equation*}
\leq \left\vert V\right\vert \int f^{p}\left( X\left( \tau \right)
,v^{\prime },\mathbf{Y}\left( \tau \right) ,\tau \right) \mbox{d}v^{\prime }%
\mbox{d}X\mbox{d}\mathbf{Y}.
\end{equation*}%
Taking the $p$-th power of (\ref{ineqpom}) and integrating over $x,$ $v,$
and $\mathbf{y}$ yields 
\begin{equation}
\Vert f(t)\Vert \hbox{\raise -1mm \hbox{$_{L^{p}}$}}\leq   \label{pom1}
\end{equation}%
\begin{equation*}
\leq \Vert f_{0}\Vert \hbox{\raise -1mm \hbox{$_{L^{p}}$}}+C\left\{ 1+\Psi
\left( C\left[ 1+\Phi \left( \sup_{0\leq s\leq t}\big\vert\mathbf{S}\left(
X\left( s\right) ,s\right) \big\vert\right) \right] \right) \right\} \times
\!\int_{0}^{t}\Vert f\left( \tau \right) \Vert 
\hbox{\raise -1mm
\hbox{$_{L^{p}}$}}\mbox{d}\tau .
\end{equation*}%
Lemma \ref{S-gradS} implies 
\begin{equation}
\sup_{0\leq \tau \leq t}\Vert \mathbf{S}\left( \mathbf{\cdot },\tau \right)
\Vert \hbox{\raise -1mm \hbox{$_{L^{\infty}}$}}\leq C\sup_{0\leq \tau \leq
t}\Vert n\left( t\right) \Vert \hbox{\raise -1mm
\hbox{$_{L^1}$}}=C\Vert n(0)\Vert \hbox{\raise
-1mm \hbox{$_{L^1}$}}\leq C\Vert f_{0}\Vert 
\hbox{\raise -1mm
\hbox{$_{L^1}$}}.  \label{pom2}
\end{equation}%
Consequently, using (\ref{pom1}) and (\ref{pom2}), we obtain 
\begin{equation}
\Vert f(t)\Vert \hbox{\raise -1mm \hbox{$_{L^p}$}}\leq \Vert f_{0}\Vert %
\hbox{\raise -1mm \hbox{$_{L^{p}}$}}+C\Big\{1+\Psi \big(C[1+\Phi (C\Vert
f_{0}\Vert \hbox{\raise -1mm \hbox{$_{L^1}$}})]\big)\Big\}\times
\int_{0}^{t}\Vert f(\tau )\Vert \hbox{\raise -1mm \hbox{$_{L^p}$}}\mbox{d}%
\tau .  \label{pom3}
\end{equation}%
Applying the Gronwall inequality, we obtain, for all $1\leq p\leq \infty ,$ 
\begin{equation}
\Vert f(t)\Vert \hbox{\raise -1mm \hbox{$_{L^p}$}}\leq C\left(
k_{i},k_{i}^{0},d_{i},t_{e},t_{a},\Vert f_{0}\Vert 
\hbox{\raise -1mm
\hbox{$_{L^{p}}$}},\mathop{\mbox{supp}}f_{0},\Psi ,\Phi ,|V|\right) <\infty .
\label{pom4}
\end{equation}%
We now compute a priori estimates on derivatives of $f$. We differentiate (%
\ref{hyperpsim}) with respect of $x$, integrate along the characteristic (%
\ref{X}) -- (\ref{Y}) from $0$ to $t$ and use (\ref{growthassumptions}) to
get%
\begin{equation*}
\left\vert \frac{\partial f}{\partial x}\left( x,v,\mathbf{y},t\right)
\right\vert \;\;\leq \;\left\vert \frac{\partial f_{0}}{\partial x}\left(
X\left( 0\right) ,v,\mathbf{Y}\left( 0\right) \right) \right\vert
+\;C\int_{0}^{t}\Psi \left( \left\vert \mathbf{Y}(\tau )\right\vert \right)
\times 
\end{equation*}%
\begin{equation*}
\times \left[ |V|\left\vert \frac{\partial f}{\partial x}\left( X(\tau ),v,%
\mathbf{Y}(\tau ),\tau \right) \right\vert +\int_{V}\left\vert \frac{%
\partial f}{\partial x}\left( X(\tau ),v^{\prime },\mathbf{Y}(\tau ),\tau
\right) \right\vert \mbox{d}v^{\prime }\right] \mbox{d}\tau +
\end{equation*}%
\begin{equation*}
+\int_{0}^{t}\big\vert\nabla _{\mathbf{y}}\cdot \mathbf{F}\left( \mathbf{S(}%
X\left( \tau \right) ),\mathbf{Y}\left( \tau \right) \right) \big\vert%
\left\vert \frac{\partial f}{\partial x}\left( X\left( \tau \right) ,v,%
\mathbf{Y}\left( \tau \right) ,\tau \right) \right\vert \mbox{d}\tau +
\end{equation*}%
\begin{equation*}
+\int_{0}^{t}\left\vert \frac{\partial g}{\partial \mathbf{z}}\left( \mathbf{%
S}(X\left( \tau \right) )\right) \right\vert \left\vert \frac{\partial 
\mathbf{S}}{\partial x}\left( X\left( \tau \right) ,\tau \right) \right\vert
\left\vert \nabla _{\mathbf{y}}f\left( X\left( \tau \right) ,v,\mathbf{Y}%
\left( \tau \right) ,\tau \right) \right\vert \mbox{d}\tau .
\end{equation*}%
Similarly, differentiating (\ref{hyperpsim}) with respect of $y_{1}$ or $%
y_{2}$, integrating along the characteristic (\ref{X}) -- (\ref{Y}) from $0$
to $t$ and using (\ref{growthassumptions}), we obtain 
\begin{equation*}
\left\vert \nabla_{\mathbf{y}}f\left( x,v,\mathbf{y},t\right) \right\vert
\;\;\leq \;\left\vert \nabla _{\mathbf{y}}f_{0}\left( X(0),v,\mathbf{Y}%
(0)\right) \right\vert +C\int_{0}^{t}\Psi \left( \left\vert \mathbf{Y(}%
X\left( \tau \right) ,\tau )\right\vert \right) \times \;
\end{equation*}%
\begin{equation*}
\times \Big\{|V|\big[\left\vert f\right\vert +\left\vert \nabla _{\mathbf{y}%
}f\right\vert \big]\left( X(\tau ),v,\mathbf{Y}(\tau ),\tau \right) +\int_{V}%
\big[\left\vert f\right\vert +\left\vert \nabla _{\mathbf{y}}f\right\vert %
\big]\left( X(\tau ),v^{\prime },\mathbf{Y}(\tau ),\tau \right) dv^{\prime }%
\Big\}\mbox{d}\tau +
\end{equation*}%
\begin{equation*}
+C\int_{0}^{t}\big\vert\nabla _{\mathbf{y}}\cdot \mathbf{F}\left( \mathbf{S(}%
X\left( \tau \right) ),\mathbf{Y}\left( \tau \right) \right) \big\vert%
\big\vert\nabla _{\mathbf{y}}f\left( X\left( \tau \right) ,v,\mathbf{Y}%
\left( \tau \right) ,\tau \right) \big\vert\mbox{d}\tau .
\end{equation*}%
If the interior of set $V$ is nonempty, we can also define the derivatives
of $f$ with respect of $v$ for any point in the interior of set $V.$
Differentiating (\ref{hyperpsim}) with respect of $v$ and integrating along
the characteristic (\ref{X}) -- (\ref{Y}) from $0$ to $t$, it implies 
\begin{equation*}
\left\vert \frac{\partial f}{\partial v}\left( x,v,\mathbf{y},t\right)
\right\vert \;\;\leq \;\left\vert \frac{\partial f_{0}}{\partial v}\left(
X\left( 0\right) ,v,\mathbf{Y}\left( 0\right) \right) \right\vert
+\int_{0}^{t}\left\vert \frac{\partial f}{\partial x}\left( X\left( \tau
\right) ,v,\mathbf{Y}\left( \tau \right) ,\tau \right) \right\vert \mbox{d}%
\tau +
\end{equation*}%
\begin{equation*}
+\;C\int_{0}^{t}\Psi \left( \left\vert \mathbf{Y}(\tau )\right\vert \right)
\times 
\bigg\{
|V|
\left\vert 
f\left( X(\tau ),v,\mathbf{Y}(\tau),\tau \right)
\right\vert
+
\end{equation*}%
\begin{equation*}
+
|V|
\left\vert 
\frac{\partial f}{\partial v}
\left( X(\tau ),v,\mathbf{Y}(\tau ),\tau \right)
\right\vert 
+
\int_{V}
\left\vert 
f
\left( X(\tau ),v^{\prime },\mathbf{Y}(\tau),\tau \right) 
\right\vert 
dv^{\prime }
\bigg\}\mbox{d}\tau 
+
\end{equation*}%
\begin{equation*}
+\int_{0}^{t}\big\vert\nabla _{\mathbf{y}}\cdot \mathbf{F}\left( \mathbf{S(}%
X\left( \tau \right) ),\mathbf{Y}\left( \tau \right) \right) \big\vert%
\left\vert \frac{\partial f}{\partial v}\left( X\left( \tau \right) ,v,%
\mathbf{Y}\left( \tau \right) ,\tau \right) \right\vert \mbox{d}\tau .
\end{equation*}%
Using (\ref{pom4}), (\ref{det}), Lemma \ref{S-gradS} and Gronwall
inequality, we deduce%
\begin{equation}
\left\Vert \frac{\partial f}{\partial x}(t)\right\Vert _{L^{p}}+\left\Vert 
\frac{\partial f}{\partial v}(t)\right\Vert _{L^{p}}+
\left\Vert 
\nabla_{\mathbf{y}}f(t)
\right\Vert _{L^{p}}\leq 
\label{derestimates}
\end{equation}%
\begin{equation*}
\leq C\left( k_{i},k_{i}^{0},d_{i},t_{e},t_{a},\Vert f_{0}\Vert 
\hbox{\raise
-1mm \hbox{$_{W^{1,p}}$}},\Vert \mathbf{S}\left( 0\right) \Vert 
\hbox{\raise
-1mm \hbox{$_{W^{1,p}}$}},\mathop{\mbox{supp}}f_{0},\mathop{\mbox{supp}}%
\mathbf{S}_{0},\Psi ,\Phi ,|V|\right) <\infty .
\end{equation*}%
Combining (\ref{pom4}) and (\ref{derestimates}), we obtain (\ref{spacef}).
Using Lemma \ref{S-gradS}, we get the estimate (\ref{spaceS}).

\rightline{Q.E.D}

\medskip

\noindent \textbf{Remark.} Using Sobolev embedding theorems, we get global
existence of classical solutions provided that initial data are smooth.

\section{Biological background}

\label{secbiodetails}

In order to study the general system (\ref{hyperp}) -- (\ref{parabS}), we
have to first specify realistic assumptions on the parameters of the model.
To this end, we summarize the relevant biological processes in Section \ref%
{secbache} and we extract the mathematical assumptions in Section \ref%
{secmajorassumptions}. These assumptions will be later used to prove the
global existence results in Section \ref{secglobalexistence}.

\subsection{Bacterial chemotaxis}

\label{secbache}

As discussed before, the bacterial movement can be viewed as a biased random
walk. Bacterial motility is commonly provided by flagella, which are long,
spiral-shaped protein rods that stick out from the surface of the cell \cite%
{Salyers:2001:MDD}. The example of flagellated bacterium is the enteric
bacterium \textit{E.coli} which has 6-8 flagella. It has two modes of
behaviour based on counterclockwise and clockwise flagellar rotation. When
the flagella rotate counterclockwise (CCW), they all point in one direction
and consequently the cell moves forward in a straight \textquotedblleft
run". The speed of running is $s=10-20\mu \hbox{m/sec}.$ Clockwise (CW)
rotation of the flagella causes the flagella to point in different
directions, and the cell tumbles in place. Tumbling reorients the cell, so
that it can move in new direction when running starts again.

\begin{figure}[tbp]
\centerline{\epsfxsize=4.5in\epsfbox{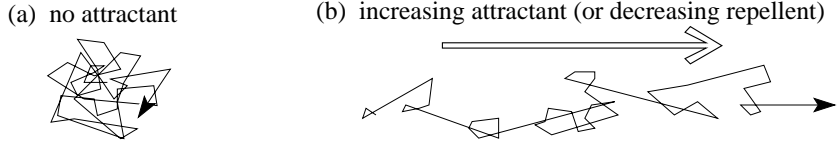}}
\caption{(a) \textit{A typical bacterial trajectory when no attractant is
present.} (b) \textit{Under the influence of an attractant, the cell
increases its time in running in a favourable direction.}}
\label{figrw}
\end{figure}

For \textit{E.coli}, the duration of both runs and tumbles are exponentially
distributed with means of 1 sec and $10^{-1}$ sec respectively if an
extracellular chemical signal is not present \cite{Block:1983:AKB}. Under
the influence of an attractant, the cell increases its time in running in a
favourable direction -- see Figure \ref{figrw}. 
As the mean time for tumbling is ten times smaller than the mean time of
running, we can often neglect the time spent tumbling and we can model the
movement of the bacterium as a velocity jump process \cite%
{Othmer:1988:MDB,Erban:2004:ICB,Erban:2004:STS} as we already did in Section %
\ref{secintro}. It means that the bacterium runs in some direction and at
random instants of time it changes its direction with mean turning rate $%
\lambda(\mathbf{y}).$

Since the bacteria move with more or less constant speed, the set $V$ of all
available velocities might be considered equal to $V=s \mathcal{S}^{N-1}$
where $\mathcal{S}^{N-1}$ is a unit sphere in $\mathbb{R}^{N}$ and $s$ is
the speed of the bacterium. Let us note that set $V=s \mathcal{S}^{N-1}$
satisfies the general condition (\ref{setvel}) (the presented theory works
for any set $V$ which satisfy (\ref{setvel})).

The kernel $K(\mathbf{v},\mathbf{v}^{\prime },\mathbf{y})$ gives the
probability of a change in velocity from $\mathbf{v}^{\prime}$ to $\mathbf{v}
$, given that a reorientation occurs. The simples possibility is to assume
that kernel is constant, i.e. 
\begin{equation}
K(\mathbf{v},\mathbf{v}^{\prime },\mathbf{y})=\frac{1}{|V|}.
\label{romundistr1}
\end{equation}%
This formula satisfies the normalization condition (\ref{normalcond}). The
underlying assumption behind (\ref{romundistr1}) is that (during the tumble)
bacterium simply choose a new direction randomly which is relatively a good
approximation for the bacterial chemotaxis, although there is also some bias
in the direction of the preceding run \cite{Berg:1972:CEC,Berg:1975:HBS}.
More realistically, one can assume that the turning kernel is a function of
the angle between new and old velocity, i.e. 
\begin{equation}
K(\mathbf{v},\mathbf{v}^{\prime },\mathbf{y})= k(\theta), \qquad \mbox{where}
\; \; \cos(\theta) = \frac{\mathbf{v} \cdot \mathbf{v^{\prime}}}{|\mathbf{v}%
| \, |\mathbf{v^{\prime}}|}.  \label{romundistr2}
\end{equation}%
Whatever the choice of $K(\mathbf{v},\mathbf{v}^{\prime },\mathbf{y})$ is,
we may assume that it is bounded from above by a constant, i.e. 
\begin{equation}
K(\mathbf{v},\mathbf{v}^{\prime},\mathbf{y})\le C  \label{romundistr3}
\end{equation}%
where $C$ is independent of $\mathbf{v},$ $\mathbf{v}^{\prime}$ and $\mathbf{%
y}$. Next, we have to specify the choice of (\ref{rom14}) and the properties
of the turning frequency $\lambda(\mathbf{y}).$

Chemotaxis is the process by which a cell alters its movement in response to
an extracellular chemical signal. From the microscopic (cell) point of view,
bacterial chemotaxis consists of several steps. First, the cell detects the
signal using its receptors. Then the signal information propagates through
the signal transduction biochemical network described by (\ref{rom14}). The
output of this network is a phosphorylated form of the protein CheY (denoted
CheY-P) which alters the motor behaviour of the flagellar motors and
consequently, the movement of the cell. CCW is the default state in the
absence of CheY-P, which binds to motor proteins and increases CW rotation.
Attractant binding to a receptor reduces the phosphorylation rate of CheY
and thereby increases the time spent in running state which constitutes the
fast response to a signal called \textit{excitation} of signal transduction
network. Another important aspect of signal transduction network is \textit{%
adaptation} which means that the response (probability per unit time of
CCW/CW rotation of flagella) returns to baseline levels on a time scale that
is slow compared to excitation, provided that there is no further change in
attractant concentration around the cell.

A schematic of the signal transduction pathway is shown in Figure \ref%
{pathways} and it can be described as follows \cite%
{Spiro:1997:MEA,Stock:1996:C,Erban:2004:STS}. 
\begin{figure}[th]
\centerline{\epsfxsize=5in\epsfbox{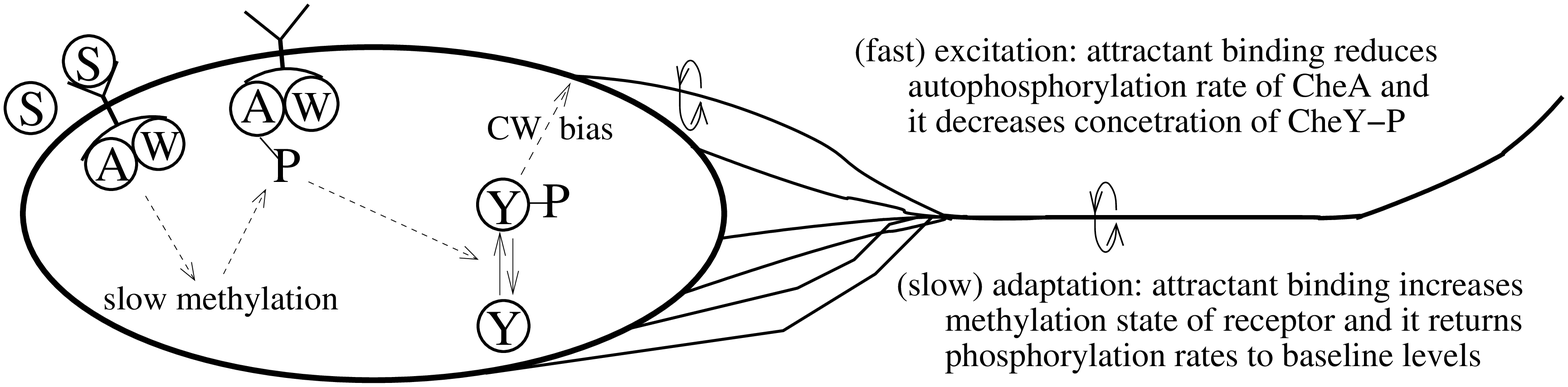}}
\caption{\textit{Excitation and adaptation in signal transduction pathway of
E. coli (from \protect\cite{Erban:2005:ICB}, with permission).}}
\label{pathways}
\end{figure}
Aspartate, the attractant most commonly-used in experiments (denoted S in
Figure \ref{pathways}), binds directly to the periplasmic domain of its
receptor, Tar. The cytoplasmic domain of Tar forms a stable complex with the
signaling proteins CheA and CheW (denoted A and W, respectively, in Figure %
\ref{pathways}), and the stability of this complex is not affected by ligand
binding \cite{Gegner:1992:AMR}. The signaling currency is in the form of
phosphoryl groups (-P), made available to the CheY (denoted Y in Figure \ref%
{pathways}) and CheB (not shown in Figure \ref{pathways}) through
autophosphorylation of CheA. Receptor complexes have two alternative
signaling states. In the attractant-bound form, the receptor inhibits CheA
autokinase activity; in the unliganded form, the receptor stimulates CheA
activity. Consequently, the response of the signal transduction network to a
step increase of the attractant concentration is as follows. First, the
attractant binding to a receptor reduces the autophosphorylation rate of
CheA. The level of phosphorylated CheA is thus lowered, causing less
phosphate to be transferred to CheY, yielding a lowered level of CheY-P. As
a result, tumbling is suppressed, and the cell's run length increases. This
constitutes the \textit{excitation} response of the system. Next slow
methylation and demethylation reactions begin to influence the response.
Ligand-bound receptors are more readily methylated than unliganded
receptors, and the lowered level of CheA-P causes a decrease in the level of
CheB-P, thereby reducing its demethylation activity. As a result, the
equilibrium of the system shifts in the direction of the higher methylation
states. The autophosphorylation rate of CheA is faster when the associated
Tar-CheA-CheW complex is in a higher methylation state, and so there is
finally a shift back toward the receptor states containing CheA-P. As a
result, CheY-P returns to its prestimulus level, and thus so does the CW
bias of the cell. This constitutes the \textit{adaptation} response. These
key steps, excitation via reduction in CheY-P, when a receptor is occupied,
and adaptation via methylation of the receptors, have been already
incorporated in the mathematical models of the bacterial signal transduction 
\cite{Spiro:1997:MEA,Barkai:1997:RSB,Morton-Firth:1999:FEB}.

Since the turning rate of bacterium is altered by CheY \cite{Cluzel:2000:UBM}%
, we can write $\lambda (\mathbf{y})\equiv \lambda (y_{1})$ where $y_{1}$
denotes the concentration of the phosphorylated form of CheY. Hence, the
individual-based model for bacterial chemotaxis is fully specified by the
equation (\ref{rom14}) which is integrated along the trajectory of each
cell, and by the $y_{1}$ component of the solution together with $\lambda
=\lambda (y_{1})$. The essential aspects of the dynamics which must be
captured by model (\ref{rom14}) are (i) it must exhibit excitation, which
here means a change in the turning frequency $\lambda (y_{1})$ in response
to a stimulus, (ii) the bias must return to baseline levels (i.e., the
response must adapt) on a time scale that is slow compared to excitation,
and (iii) the signal transduction network should amplify signals
appropriately \cite{Bourret:1991:STP,Segall:1986:TCB}. The mathematical
assumptions on (\ref{rom14}) and $\lambda (y_{1})$ are given in Section \ref%
{secmajorassumptions}.

\subsection{Mathematical assumptions on the signal transduction network}

\label{secmajorassumptions}

The mathematical model of the signal transduction network (\ref{rom14}) 
can be rewritten in the following form 
\begin{equation}
\frac{\mbox{d}\mathbf{y}}{\mbox{d}t}=\mathbf{F}(\mathbf{C}(t),\mathbf{y}%
)\qquad \mbox{where}\qquad \mathbf{C}(t)=\mathbf{S}(\mathbf{x}(t),t).
\label{rom14new}
\end{equation}%
The vector function $\mathbf{C}(t)$ gives signal values which are seen by a
cell along its trajectory. Time evolution of $\mathbf{y}$ in equation (\ref%
{rom14new}) is controlled by the input time dependent vector $\mathbf{C}(t).$
Therefore, it is natural to describe the behaviour of $\mathbf{F}$ in terms
of the input function $\mathbf{C}(t).$

The mathematical formulation of the \textit{adaptation} property of the
signal transduction network (\ref{rom14}) can be written in the following
form. There exists a universal constant $\overline{y}_{1}$ such that for any
constant signal along the trajectory $\mathbf{C}_{0}$, i.e. $\mathbf{C}%
(t)\equiv \mathbf{C}_{0}=\mbox{{\rm const}}$ and for any initial condition $%
\mathbf{y}(0)=\mathbf{y}_{0}$, the solution of the system $(\ref{rom14})$
satisfies 
\begin{equation}
\lim_{t\rightarrow \infty }y_{1}(t)=\overline{y}_{1}.  \label{adaptproperty}
\end{equation}%
Formula (\ref{adaptproperty}) describes the perfect adaptation. From the
application point of view, it is desirable that the signal transduction
model satisfies (at least approximately) the adaptation property for a
reasonably large set of signals. However, the existence theorems presented
in Section \ref{secglobalexistence} do not require perfect adaptation and we
will prove the existence of solutions even for models which do not satisfy (%
\ref{adaptproperty}). It is worthwhile to note that the simplified model of
excitation-adaptation dynamics (\ref{simmodel}) from Section \ref{secmotiv}
satisfied adaptation property (\ref{adaptproperty}). In fact, $%
y_{1}(t)\rightarrow 0$ as $t\rightarrow \infty $ for any constant signal,
i.e. $y_{1}$ adapts perfectly to any constant stimulus. Moreover, model (\ref%
{simmodel}) describes the excitation-adaptation dynamics as discussed in
Section \ref{secbache} provided that we choose $t_{e}<t_{a}$. Here, the time
constants $t_{e}$ and $t_{a}$ are labeled in anticipation of using $y_{1}$
for the internal response, and $y_{2}$ as the adaptation variable, and
therefore we call $t_{e}$ and $t_{a}$ the excitation and adaptation time
constant, respectively \cite{Erban:2004:ICB}.

\begin{figure}[th]
\centerline{\epsfxsize=4.8in\epsfbox{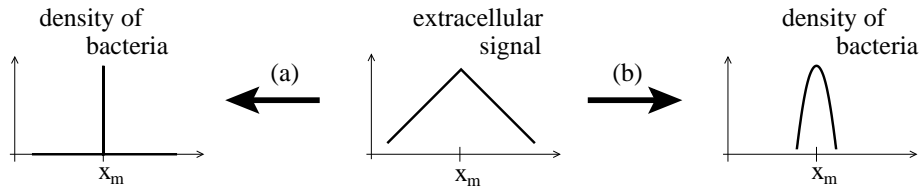}}
\caption{\textit{Schematic of behaviour of hypothetical cells which
\textquotedblleft perfectly avoid going in wrong directions" (panel on the
left) and hypothetical cells which \textquotedblleft perfectly follow good
directions" (panel on the right). Details are explained in the text. }}
\label{figexpblow}
\end{figure}

In order to model the random walk of the individual bacterium, we must have
a good understanding of the dependence of the (output) turning rate $\lambda
(y_{1})$ on the (input) signal function $\mathbf{C}(t).$ If the input signal
function is constant then the behaviour of $\lambda (y_{1})$ follows the
adaptation property. On the other hand, time dependent input $\mathbf{C}(t)$
can introduce large variations in $\lambda (y_{1}).$ The time derivative of $%
\mathbf{C}(t),$ i.e. the time derivative of the signal seen by a cell, is
equal to 
\begin{equation}
\frac{\mbox{d}\mathbf{C}}{\mbox{d}t}=\mathbf{v}\cdot \frac{\partial \mathbf{S%
}}{\partial \mathbf{x}}+\frac{\partial \mathbf{S}}{\partial t}.
\label{timdercell}
\end{equation}%
To see what type of conditions on the turning rate $\lambda$
are reasonable, let us consider the time
independent signal (attractant) with a maximum at the point $x_{m}$ as it is
schematically shown in one dimension in Figure \ref{figexpblow} (panel in
the middle). 
We consider that bacteria move with the fixed speed either to the right or
left and we discuss the following two simple cases of dependence of output $%
\lambda (y_{1})$ on input $\mathbf{C}(t)$. 
\begin{equation*}
\mbox{(a)}\;\;\lambda (y_{1})=\left\{ 
\begin{array}{ll}
1\quad  & \mbox{for}\;\mbox{d}C/\mbox{d}t\geq 0; \\ 
\infty \quad  & \mbox{for}\;\mbox{d}C/\mbox{d}t<0;%
\end{array}%
\right. \qquad \mbox{(b)}\;\;\lambda (y_{1})=\left\{ 
\begin{array}{ll}
0\quad  & \mbox{for}\;\mbox{d}C/\mbox{d}t\geq 0; \\ 
1\quad  & \mbox{for}\;\mbox{d}C/\mbox{d}t<0.%
\end{array}%
\right. 
\end{equation*}%
Let us note that cases (a) and (b) are considered as definitions of the
input-output behaviour in two extreme cases (these definitions are not
connected with any underlying differential equation in this example).

First, suppose that a bacterium is at the position $x<x_{m}$. If we use
input-output behaviour (a), then the cell goes to the right. It sometimes
\textquotedblleft turns" to the left but it instantly turns back. So, the
cell spends all the time going to the right, and case (a) is an example of
the individual-based model where cells perfectly avoid going in wrong
directions. If we use input-output behaviour (b), then the right going cells
never turn (for $x<x_{m}$). Hence, case (b) is an example of the 
individual-based model
where cells perfectly follow good directions. Both cases (a) and (b)
describe the simple transport of bacteria for $x<x_{m}$. The difference of
these models is when cells reach the maximum of the signal $x_{m}$. In case
(a), cells instantly turn back. It means that the final positions of all
bacteria are equal to $x_{m}$ and a Dirac-like distribution is created in
finite time (see Figure \ref{figexpblow}, panel on the left). In case (b),
cells continue movement to the region $x>x_{m}$ and the final distribution
profile is smooth, as shown schematically in Figure \ref{figexpblow} (panel
on the right).

The previous simple example shows that singularities might develop if the
turning rate is too large (without a reasonable control by the signal
change), as in case (a) where cells perfectly avoid going in wrong
directions. This observation suggests for growth conditions on $\lambda
(y_{1})$ from above which prevent formation of singularities. The necessary
conditions on the turning frequency $\lambda (y_{1})$ is $\lambda
(y_{1})\geq 0$ and our heuristic conclusions can be incorporated to the
following growth estimate 
\begin{equation}
\lambda (y_{1})\leq C\left( 1+\Lambda \left( \left\vert \mathbf{C}%
\right\vert \right) +\left\vert \frac{\mbox{d}\mathbf{C}}{\mbox{d}t}%
\right\vert \right) ,  \label{growthlambda}
\end{equation}%
where $\Lambda \left( \cdot \right) \in C\left( \mathbb{R}\right) $ is a
non-negative, nondecreasing continuous function. The verification of growth
estimate (\ref{growthlambda}) depends on the particular form of $\mathbf{F}%
(\cdot ,\cdot )$ and $\lambda (\cdot ).$ For example, if (\ref{rom14}) and $%
\lambda (\cdot )$ satisfy 
\begin{equation}
|y_{1}|\leq C_{1}\left( 1+\left\vert \frac{\mbox{d}\mathbf{C}}{\mbox{d}t}%
\right\vert ^{\omega }\right) ,\qquad \lambda (y_{1})\leq
C_{2}(1+|y_{1}|^{\sigma }),\qquad \omega \sigma \leq 1,  \label{growthylam}
\end{equation}%
then (\ref{growthlambda}) follows. There are several other conditions on $%
\mathbf{F}(\cdot ,\cdot )$ and $\lambda (\cdot )$ which also guarantee
growth estimate (\ref{growthlambda}). Hence, we do not formulate our growth
estimates in terms of $\mathbf{F}(\cdot ,\cdot )$ and $\lambda (\cdot )$,
but we simply assume (\ref{growthlambda}) directly in our existence
theorems. Using formula (\ref{kernelT}), we can formulate the estimate (\ref%
{growthlambda}) also in terms of the kernel $T(\mathbf{v},\mathbf{v}^{\prime
},\mathbf{y}).$

Using estimate (\ref{romundistr3}) and definition (\ref{kernelT}), we can
write the growth assumption on $T$ in the following form 
\begin{equation}
T(v,v^{\prime },\mathbf{y}) \leq C |\lambda (y_{1})|.  \label{growthT}
\end{equation}%
We also have to assume a growth assumption of ${\nabla }_{\mathbf{y}}{\cdot }%
\mathbf{F}$. In Theorem \ref{theorem2}, we assume that there exists a
non-negative, nondecreasing continuous function $\Pi \left( \cdot \right)
\in C\left( \mathbb{R}\right) $ satisfying%
\begin{equation}
\left\vert {\nabla }_{\mathbf{y}}{\cdot }\mathbf{F}\left( \mathbf{z,y}%
\right) \right\vert \leq C\left( 1+\Pi \left( \left\vert \mathbf{z}%
\right\vert \right) \right).  \label{divF2}
\end{equation}
Notice that our simple model (\ref{simmodel}) satisfies (\ref{divF2}). A
different condition on ${\nabla }_{\mathbf{y}}{\cdot }\mathbf{F}$ is studied
also in Corollary \ref{cor1}.

\subsection{Mathematical assumptions on the dynamics of the extracellular
signals}

\label{secmajorassumptionssignal}

Various forms of $\mathbf{R}(\mathbf{S},n)$ can be considered. The simplest
case from the mathematical point of view is when the extracellular signals
are nutrients which are consumed by cells, i.e. 
\begin{equation}
\mathbf{R}(\mathbf{S},n)= - \mathbf{K} \mathbf{S} n  \label{eatonly}
\end{equation}
where $\mathbf{K}$ is a diagonal nonnegative $M \times M$ matrix (with rate
constants on the diagonal). One can also assume that the cells produce
signals which are degraded at some rate, i.e. 
\begin{equation}
\mathbf{R}(\mathbf{S},n)= n [k_1,k_{2},\dots ,k_{M}]^{T} - \mathbf{K} 
\mathbf{S}  \label{proddegr}
\end{equation}
where $k_1,$ \dots, $k_M$ are rates of production of the different
components of the signal and $\mathbf{K}$ is a diagonal nonnegative $M
\times M$ matrix. If we allow the nondiagonal terms in matrix $\mathbf{K}$,
then the extracellular coupling of the signals (e.g. reactions between signals)
is added to the model. One can also consider that some signals can be
produced by cells and some signals can be degraded by cells, i.e.
effectively combining (\ref{eatonly}) and (\ref{proddegr}). Moreover, we can
also assume that some signals can be attractants while other signals can be
repellents etc.

Depending on the model system, there are many possibilities to specify the
dynamics of the extracellular signal. In what follows, we use (\ref{proddegr}%
). However, it is possible to modify and prove the following existence
theorems using different evolution equations for the extracellular signal
too. The only requirement is that the evolution equation for the
extracellular signal must satisfy suitable growth estimates similar to the
estimates which are proven in Lemma \ref{gradS} for (\ref{proddegr}).

\section{Global existence for the general signal transduction models}

\label{secglobalexistence}

In this section, we prove global existence results using the framework of
Sections \ref{secmajorassumptions} and \ref{secmajorassumptionssignal}. We
will work in one-dimensional physical space, i.e. $N=1$ and we first assume
the case of elliptic equations for the extracellular signals. Hence, system
of equations (\ref{hyperp}) -- (\ref{parabS}) reads as follows 
\begin{equation}
\frac{\partial f}{\partial t}+\nabla_{\!x}\cdot \mathbf{v} f + \nabla _{\!%
\mathbf{y}} \cdot \mathbf{F}(\mathbf{S}(\mathbf{x}),\mathbf{y}) f = \int_{V}
T(v,v^{\prime},\mathbf{y}) \Big[ f(v^{\prime}) - f(v) \Big] \mbox{d}%
v^{\prime},  \label{hyperpell}
\end{equation}
\begin{equation}
d_{i}\frac{\partial ^{2}S_{i}}{\partial x^{2}}+k_{i}n-k_{i}^{0}S_{i} = 0,
\qquad \;i=1,\dots ,M,  \label{ellipticS}
\end{equation}
where $d_{i}$, $k_{i}$ and $k_{i}^{0}$ are positive constants and $n\equiv
n(x,t)$ is the macroscopic density of individuals at point $x \in \mathbb{R}$
and time $t$ given by (\ref{equationforn}). Now, we can formulate the
existence theorem.

\begin{theorem}
\label{theorem2} Let us assume $(\ref{growthlambda})$, $(\ref{growthT})$ and 
$(\ref{divF2})$. Assume that $f_{0}\in L^{1}\cap L^{\infty }({}\mathbb{R}%
\times V\times \mathbb{R}^{m})$ and let initial condition $S_{0}\in \lbrack
W^{2,p}({}\mathbb{R})]^{M}$ satisfies $(\ref{ellipticS})$. Then there exists
a global solution of system $(\ref{hyperpell})$ -- $(\ref{ellipticS})$
satisfying, for all $t\geq 0$ 
\begin{equation}
f(\cdot ,\cdot ,\cdot ,t)\in L^{1}\cap L^{\infty }({}\mathbb{R}\times
V\times \mathbb{R}^{m}),  \label{fpestimate}
\end{equation}%
\begin{equation}
\mathbf{S}(\cdot ,t)\in \lbrack W^{2,p}(\mathbb{R})]^{M},\ \ \ \ \ \text{for
all }1\leq p<+\infty,  \label{Spestimate}
\end{equation}%
and initial conditions $f(\cdot ,\cdot ,\cdot ,0)=f_{0}(\cdot ,\cdot ,\cdot
) $ and $\mathbf{S}(\cdot ,0)=\mathbf{S}_{0}(\cdot )$.
\end{theorem}

\medskip

\noindent \textbf{Remark.} To avoid technicalities, we focus in Theorem \ref%
{theorem2} only on $L^{p}$ estimates of $f$. The results could be extended
to $W^{k,p}$ estimates under suitable growth assumptions on derivatives of $%
T\left( v,v^{\prime },y\right)$ and $\mathbf{F}$.

\medskip

\noindent In order to prove Theorem \ref{theorem2}, we formulate some
auxiliary lemmas. We start with the generalization of the Gronwall
inequality.

\begin{lemma}
\label{gronwall}Let $a\left( s\right) $ and $b\left( s\right) $ be positive
integrable functions on $[0,t]$. Let $w\left( t\right) $ be positive and
differentiable in $t,$ and satisfy 
\begin{equation*}
w^{\prime }\leq a\left( t\right) w\ln w+b\left( t\right) w.
\end{equation*}%
Then 
\begin{equation*}
w\left( t\right) \leq \left[ w\left( 0\right) \exp \left(
\int_{0}^{t}b\left( s\right) e^{-\int_{0}^{s}a\left( \tau \right) d\tau
}ds\right) \right] ^{\exp \left( \int_{0}^{t}a\left( s\right) ds\right) }.
\end{equation*}
\end{lemma}

\noindent \textbf{Proof. } See \cite[Lemma 4]{Hwang:2005:GEC}.

\rightline{Q.E.D.}

\noindent The characteristics of the hyperbolic equation (\ref{hyperpell})
are given for $N=1$ as the solution of (\ref{chareqn}). The back-in-time
characteristics starting at $\left( x,v,\mathbf{y},t\right) $ are given as 
\begin{eqnarray}
X\left( s;x,v\mathbf{,y},t\right)  &=&x-v\left( t-s\right) ,\ \   \label{X2}
\\
\mathbf{Y}\left( s;x,v\mathbf{,y},t\right)  &=&\mathbf{y}-\int_{s}^{t}%
\mathbf{F}\left( \mathbf{S}\left( X(\tau ),\tau \right) ,\mathbf{Y}(\tau
)\right) \mbox{d}\tau .  \label{Y2}
\end{eqnarray}%
The generalization of Lemma \ref{lemma3} is given as the following Lemma.

\begin{lemma}
\label{detGeneral} Derivation of the characteristics $(\ref{X2})$ and $(\ref%
{Y2})$ with respect to the initial conditions gives 
\begin{equation}
\frac{\partial X}{\partial x}=1 \qquad \text{and} \qquad \frac{\partial 
\mathbf{Y}}{\partial \mathbf{y}} = \exp \left[ -\int_{s}^{t}\frac{\partial 
\mathbf{F}}{\partial \mathbf{y}} \left( \mathbf{S}\left( X(\tau ),\tau
\right) ,\mathbf{Y}(\tau)\right) \mbox{d}\tau \right].
\label{dercharGeneral}
\end{equation}%
Moreover, 
\begin{equation}
\det \frac{\partial \mathbf{Y}}{\partial \mathbf{y}} = \exp \left[
-\int_{s}^{t}\nabla _{\mathbf{y}}\cdot \mathbf{F} \left( \mathbf{S}\left(
X(\tau ),\tau \right) ,\mathbf{Y}(\tau)\right) \mbox{d}\tau \right].
\label{determyyGeneral}
\end{equation}
\end{lemma}

\noindent \textbf{Proof.} We differentiate (\ref{Y2}) with respect to $%
\mathbf{y}$ to get%
\begin{equation}
\frac{\partial \mathbf{Y}}{\partial \mathbf{y}}=\mathbf{I}_{m}+\int_{t}^{s}%
\frac{\partial \mathbf{F}}{\partial \mathbf{y}}\left( \mathbf{S}\left(
X(\tau ),\tau \right) ,\mathbf{Y}(\tau )\right) \frac{\partial \mathbf{Y}}{%
\partial \mathbf{y}}\left( \tau \right) \mbox{d}\tau .  \label{derY2}
\end{equation}%
where $\mathbf{I}_{m}$ is the $m\times m$ identity matrix. Let 
\begin{equation*}
\mathbf{G}\left( s\right) =\int_{t}^{s}\frac{\partial \mathbf{F}}{\partial 
\mathbf{y}}\left( \mathbf{S}\left( X(\tau ),\tau \right) ,\mathbf{Y}(\tau
)\right) \frac{\partial \mathbf{Y}}{\partial \mathbf{y}}\left( \tau \right) %
\mbox{d}\tau ,
\end{equation*}%
then we have%
\begin{equation*}
\mathbf{G}^{\prime }\left( s\right) -\mathbf{G}\left( s\right) \frac{%
\partial \mathbf{F}}{\partial \mathbf{y}}\left( \mathbf{S}\left( X(\tau
),\tau \right) ,\mathbf{Y}(\tau )\right) =\frac{\partial \mathbf{F}}{%
\partial \mathbf{y}}\left( \mathbf{S}\left( X(\tau ),\tau \right) ,\mathbf{Y}%
(\tau )\right) .
\end{equation*}%
Integrating the last equation, we obtain (\ref{dercharGeneral}). Since the
determinant of the exponential of the matrix is the exponential of the trace
of the matrix, we have 
\begin{equation*}
\det \frac{\partial \mathbf{Y}}{\partial \mathbf{y}}=\exp \left[ %
\mathop{\mbox{{\rm trace}}}\left( -\int_{s}^{t}\frac{\partial \mathbf{F}}{%
\partial \mathbf{y}}\left( \mathbf{S}\left( X(\tau ),\tau \right) ,\mathbf{Y}%
(\tau )\right) \mbox{d}\tau \right) \right] =
\end{equation*}%
\begin{equation*}
=\exp \left[ -\int_{s}^{t}\nabla _{\mathbf{y}}\cdot \mathbf{F}\left( \mathbf{%
S}\left( X(\tau ),\tau \right) ,\mathbf{Y}(\tau )\right) \mbox{d}\tau \right]
.
\end{equation*}%
Hence, we have proved (\ref{determyyGeneral}).

\rightline{Q.E.D.}

\medskip

\noindent Next, we present the growth estimates on the extracellular signal $%
\mathbf{S}$ and on its derivatives. The time and space derivatives of the
signal vector $\mathbf{S}$ are controlled by logarithm of the $L^{2}$-norm
of the cell density. Note that the analogous result was also shown in \cite[%
Lemma 4]{Hwang:2005:GEC} for the parabolic equation for the extracellular
signal. The difference between \cite[Lemma 4]{Hwang:2005:GEC} and Lemma \ref%
{gradS} is that we prove also the estimate on the time derivative as well as
the estimate on the space derivative of the signal.

\begin{lemma}
\label{gradS}If $n\in L^{\infty }([0,\infty ):L^{1}\left( \mathbb{R}\right)
\cap L^{2}\left( \mathbb{R}\right) ),$ then the solution $\mathbf{S}$ in $(%
\ref{ellipticS})$ satisfies 
\begin{equation*}
\left\Vert \mathbf{S}\left( t\right) \right\Vert _{L^{\infty }}\leq
C\left\Vert n\left( t\right) \right\Vert _{L^{1}}=C\left\Vert n\left(
0\right) \right\Vert _{L^{1}},
\end{equation*}%
\begin{equation}
\left\Vert \frac{\partial \mathbf{S}}{\partial x}\left( t\right) \right\Vert
_{L^{\infty }}\leq C\left[ 1+\left\Vert n\left( 0\right) \right\Vert
_{L^{1}}\left\{ 1+\ln \left( \left\Vert n(t) \right\Vert _{L^{2}}+1\right)
\right\} \right] ,  \label{sxderestim}
\end{equation}%
\begin{equation}
\left\Vert \frac{\partial \mathbf{S}}{\partial t}\left( t\right) \right\Vert
_{L^{\infty }}\leq C\left[ 1+\left\Vert n\left( 0\right) \right\Vert
_{L^{1}}\left\{ 1+\ln \left( \left\Vert n(t) \right\Vert _{L^{2}}+1\right)
\right\} \right] .  \label{stderestim}
\end{equation}%
where the constant $C$ depends only on $k_{i},$ $k_{i}^{0}$, $d_{i}$ and $V$.
\end{lemma}

\noindent \textbf{Proof.} Let $1\leq i\leq M.$ Taking the Fourier transform
of $(\ref{ellipticS})$ in the $x-$variable, we obtain 
\begin{equation*}
\hat{S}_{i}\left( \xi ,t\right) =\frac{k_{i}}{d_{i}}\frac{\hat{n}\left( \xi
,t\right) }{\xi ^{2}+k_{i}^{0}/d_{i}}.
\end{equation*}%
Thus we have%
\begin{eqnarray*}
\left\Vert S_{i}(t)\right\Vert _{L^{\infty }} &\leq &\left\Vert \hat{S}%
_{i}(t)\right\Vert _{L^{1}}\leq \frac{k_{i}}{d_{i}}\left\Vert \hat{n}%
(t)\right\Vert _{L^{\infty }}\int_{-\infty }^{\infty }\frac{1}{\xi
^{2}+k_{i}^{0}/d_{i}}\mbox{d}\xi  \\
&\leq &C\left( \frac{k_{i}}{d_{i}},\frac{k_{i}^{0}}{d_{i}}\right) \left\Vert
n(t)\right\Vert _{L^{1}}=C\left( \frac{k_{i}}{d_{i}},\frac{k_{i}^{0}}{d_{i}}%
\right) \left\Vert n(0)\right\Vert _{L^{1}}.
\end{eqnarray*}%
Next we estimate the $x-$derivative of the signal as follows.%
\begin{equation*}
\left\Vert \frac{\partial S_{i}}{\partial x}\left( t\right) \right\Vert
_{L^{\infty }}\leq \left\Vert \xi \hat{S}_{i}(t)\right\Vert _{L^{1}}\leq 
\frac{k_{i}}{d_{i}}\int_{-\infty }^{\infty }\frac{\left\vert \xi \right\vert
\left\vert \hat{n}\left( \xi ,t\right) \right\vert }{\xi ^{2}+k_{i}^{0}/d_{i}%
}\mbox{d}\xi =\frac{k_{i}}{d_{i}}\left\{ I_{1}+I_{2}\right\} ,
\end{equation*}%
\begin{equation*}
\mbox{where}\quad I_{1}=\int_{\left\vert \xi \right\vert \leq \Vert
n(t)\Vert _{L^{2}}^{2}}\frac{\left\vert \xi \right\vert \left\vert \hat{n}%
\left( \xi ,t\right) \right\vert }{\xi ^{2}+k_{i}^{0}/d_{i}}\mbox{d}\xi
\quad \mbox{and}\quad I_{2}=\int_{\left\vert \xi \right\vert \geq \Vert
n(t)\Vert _{L^{2}}^{2}}\frac{\left\vert \xi \right\vert \left\vert \hat{n}%
\left( \xi ,t\right) \right\vert }{\xi ^{2}+k_{i}^{0}/d_{i}}\mbox{d}\xi .
\end{equation*}%
First, we estimate the integral $I_{1}$. We obtain%
\begin{equation*}
I_{1}\leq \left\Vert \hat{n}(t)\right\Vert _{L^{\infty }}\int_{\left\vert
\xi \right\vert \leq \Vert n(t)\Vert _{L^{2}}^{2}}\frac{\left\vert \xi
\right\vert }{\xi ^{2}+k_{i}^{0}/d_{i}}\mbox{d}\xi =
\end{equation*}%
\begin{equation*}
=\left\Vert \hat{n}(t)\right\Vert _{L^{\infty }}\ln \left( \frac{\Vert
n(t)\Vert _{L^{2}}^{4}}{k_{i}^{0}/d_{i}}+1\right) \leq \left\Vert
n(t)\right\Vert _{L^{1}}\ln \left( \frac{\Vert n(t)\Vert _{L^{2}}^{4}}{%
k_{i}^{0}/d_{i}}+1\right) 
\end{equation*}%
We use H\"{o}lder's inequality with $p=q=2$ to estimate $I_{2}$ as%
\begin{eqnarray*}
I_{2} &\leq &\left\Vert \hat{n}(t)\right\Vert _{L^{2}}\left(
\int_{\left\vert \xi \right\vert \geq \Vert n(t)\Vert _{L^{2}}^{2}}\left( 
\frac{\xi }{\xi ^{2}+k_{i}^{0}/d_{i}}\right) ^{2}\mbox{d}\xi \right) ^{1/2}
\\
&\leq &\left\Vert n(t)\right\Vert _{L^{2}}\left( \int_{\left\vert \xi
\right\vert \geq \Vert n(t)\Vert _{L^{2}}^{2}}\xi ^{-2}\mbox{d}\xi \right)
^{1/2}\leq \sqrt{2}.
\end{eqnarray*}%
By combining the estimates for $I_{1}$ and $I_{2}$, we obtain (\ref%
{sxderestim}). In order to estimate the time derivative of the extracellular
signal, we take the time derivative of $(\ref{ellipticS})$ and apply the
Fourier transform in the $x-$variable to get%
\begin{equation*}
\frac{\partial \hat{S}_{i}}{\partial t}\left( \xi ,t\right) =\frac{k_{i}}{%
d_{i}}\frac{\partial \hat{n}}{\partial t}\left( \xi ,t\right) \frac{1}{\xi
^{2}+k_{i}^{0}/d_{i}}.
\end{equation*}%
By integrating $(\ref{hyperpell})$ over $v$ and $y,$ we get%
\begin{equation*}
\frac{\partial n}{\partial t}=-\frac{\partial j}{\partial x}\quad %
\mbox{where}\quad j\left( x,t\right) =\iint_{V\mathbb{\times R}^{m}}vf\left(
x,v,\mathbf{y},t\right) \mbox{d}v\mbox{d}\mathbf{y}.
\end{equation*}%
Thus we have%
\begin{equation*}
\frac{\partial \hat{S}_{i}}{\partial t}\left( \xi ,t\right) =\frac{k_{i}}{%
d_{i}}\frac{-i\xi \hat{\jmath}\left( \xi ,t\right) }{\xi ^{2}+k_{i}^{0}/d_{i}%
}.
\end{equation*}%
Then we have%
\begin{equation*}
\left\Vert \frac{\partial S_{i}}{\partial t}(t)\right\Vert _{L^{\infty
}}\leq \left\Vert \frac{\partial \hat{S}_{i}}{\partial t}(t)\right\Vert
_{L^{1}}\leq \frac{k_{i}}{d_{i}}\int_{-\infty }^{\infty }\frac{\left\vert
\xi \right\vert \left\vert \hat{\jmath}\left( \xi \right) \right\vert }{\xi
^{2}+k_{i}^{0}/d_{i}}\mbox{d}\xi .
\end{equation*}%
Notice that%
\begin{equation*}
\left\Vert \hat{\jmath}(t)\right\Vert _{L^{\infty }}\leq \left\Vert
j(t)\right\Vert _{L^{1}}\leq \iiint_{\mathbb{R\times }V\mathbb{\times R}%
^{m}}\left\vert v\right\vert f\left( x,v,\mathbf{y},t\right) \mbox{d}x%
\mbox{d}v\mbox{d}\mathbf{y}\leq 
\end{equation*}%
\begin{equation*}
\leq C\left( V\right) \left\Vert n\left( t\right) \right\Vert
_{L^{1}}=C\left( V\right) \left\Vert n(0)\right\Vert _{L^{1}},
\end{equation*}%
\begin{equation*}
\left\Vert \hat{\jmath}(t)\right\Vert _{L^{2}}=\left\Vert j(t)\right\Vert
_{L^{2}}\leq \left( \int_{\mathbb{R}}\left( \iint_{V\mathbb{\times R}%
^{m}}\left\vert v\right\vert f\left( x,v,\mathbf{y},t\right) \mbox{d}v%
\mbox{d}\mathbf{y}\right) ^{2}\mbox{d}x\right) ^{1/2}\leq 
\end{equation*}%
\begin{equation*}
\leq C\left( V\right) \left\Vert n\left( t\right) \right\Vert _{L^{2}}
\end{equation*}%
where we have used that $V$ is compact. Using similar ideas as in the proof
of estimate (\ref{sxderestim}), we prove (\ref{stderestim}).

\rightline{Q.E.D.}

\medskip

\begin{lemma}
\label{lem7} Let $\mathbf{F}$ satisfy $(\ref{divF2})$. Then the
characteristics $(\ref{X2})$ -- $(\ref{Y2})$ satisfy for all $0\leq s\leq t,$
\begin{equation}
\left[ \det \frac{\partial \mathbf{Y}}{\partial \mathbf{y}}\left( s\right) %
\right] ^{-1}\leq \exp \big[Ct\,\big\{1+\Pi \left( C\left\Vert
f_{0}\right\Vert _{L^{1}}\right) \big\}\big].  \label{deter}
\end{equation}
\end{lemma}

\noindent \textbf{Proof.} 
Using Lemma \ref{detGeneral} and (\ref{divF2}), we obtain 
\begin{equation*}
\left[ \det \frac{\partial \mathbf{Y}}{\partial \mathbf{y}}(s)\right]
^{-1}=\exp \left[ \int_{s}^{t}\nabla _{\mathbf{y}}\cdot \mathbf{F}\left( 
\mathbf{S}\left( X(\tau ),\tau \right) ,\mathbf{Y}(\tau )\right) \mbox{d}%
\tau \right] \leq
\end{equation*}%
\begin{equation*}
\leq \exp \left[ C\int_{s}^{t}1+\Pi \left( \left\vert \mathbf{S}\left(
X(\tau ),\tau \right) \right\vert \right) \mbox{d}\tau \right] \leq 
\exp \left[
C t \left\{
1+\Pi 
\left( 
\sup_{0\leq \tau \leq t}\left\Vert \mathbf{S}\left( \tau \right) 
\right\Vert _{L^{\infty }}
\right) 
\right\}
\right].
\end{equation*}%
Using Lemma \ref{gradS}, we deduce (\ref{deter}).

\vskip -2mm

\rightline{Q.E.D.}

\medskip

\noindent \textbf{Proof of Theorem 2.} Using $(\ref{growthlambda})$ and $(%
\ref{growthT})$, we obtain 
\begin{equation}
T(v,v^{\prime },\mathbf{y})\leq C\left( 1+\Lambda \left( \left\vert \mathbf{C%
}\right\vert \right) +\left\vert \frac{\mbox{d}\mathbf{C}}{\mbox{d}t}%
\right\vert \right) .  \label{estimT}
\end{equation}%
Integrating (\ref{hyperpell}) along the characteristic (\ref{X2}) -- (\ref%
{Y2}) from $0$ to $t$ and using (\ref{estimT}), we obtain 
\begin{equation*}
f\left( x,v,\mathbf{y},t\right) \leq f_{0}\left( X\left( 0\right) ,v,\mathbf{%
Y}\left( 0\right) \right) +
\end{equation*}%
\begin{equation*}
+\;C\left( V\right) \int_{0}^{t}
\bigg\{\left( 1+\left[ \Lambda \left(
\left\vert \mathbf{S}\right\vert \right) +\left\vert \frac{\partial \mathbf{S%
}}{\partial t}\right\vert +\left\vert \frac{\partial \mathbf{S}}{\partial x}%
\right\vert \right] \left( X\left( \tau \right) ,\tau \right) \right) \times
\end{equation*}%
\begin{equation*}
\times \int_{V}f\left( X(\tau ),v^{\prime },\mathbf{Y}(\tau ),\tau \right) %
\mbox{d}v^{\prime }\bigg\}\mbox{d}\tau +
\end{equation*}%
\begin{equation*}
+\int_{0}^{t}\big\vert\nabla _{\mathbf{y}}\cdot \mathbf{F}\left( \mathbf{S(}%
X\left( \tau \right) ),\mathbf{Y}\left( \tau \right) \right) \big\vert %
f\left( X\left( \tau \right) ,v,\mathbf{Y}\left( \tau \right) ,\tau \right) %
\mbox{d}\tau ,
\end{equation*}%
where we used that $V$ is compact. By virtue of assumption (\ref{divF2}), $%
\left\vert \nabla _{\mathbf{y}}\cdot \mathbf{F}\right\vert $ is bounded by 
\hfill \break 
$%
C\left( 1+\Pi \left( \left\vert \mathbf{S}\left( X(\tau ),\tau \right)
\right\vert \right) \right) .$ Thus we have 
\begin{equation}
f\left( x,v,\mathbf{y},t\right) \leq f_{0}\left( X\left( 0\right) ,v,\mathbf{%
Y}\left( 0\right) \right) +  \label{fint}
\end{equation}%
\begin{equation*}
+\;C\left( V\right) \int_{0}^{t}
\bigg\{\left( 1+\left[ \Lambda \left(
\left\vert \mathbf{S}\right\vert \right) +\left\vert \frac{\partial \mathbf{S%
}}{\partial t}\right\vert +\left\vert \frac{\partial \mathbf{S}}{\partial x}%
\right\vert \right] \left( X\left( \tau \right) ,\tau \right) \right) \times
\end{equation*}%
\begin{equation*}
\times \int_{V}f\left( X(\tau ),v^{\prime },\mathbf{Y}(\tau ),\tau \right) %
\mbox{d}v^{\prime }\bigg\}\mbox{d}\tau +
\end{equation*}%
\begin{equation*}
+C(V)\int_{0}^{t}\left( 1+\Pi \left( \left\vert \mathbf{S}\left( X(\tau
),\tau \right) \right\vert \right) \right) f\left( X\left( \tau \right) ,v,%
\mathbf{Y}\left( \tau \right) ,\tau \right) \mbox{d}\tau .
\end{equation*}%
Using Lemma \ref{lem7}, we obtain for $t\geq 0,$
\begin{equation*}
\int_{\mathbb{R}\times V\times \mathbb{R}^{M}}\int_{V}f^{p}\left( X\left(
\tau \right) ,v^{\prime },\mathbf{Y}\left( \tau \right) ,\tau \right)
dv^{\prime }\mbox{d}x\mbox{d}v\mbox{d}\mathbf{y}=
\end{equation*}%
\begin{equation*}
=\left\vert V\right\vert \int f^{p}\left( X\left( \tau \right) ,v^{\prime },%
\mathbf{Y}\left( \tau \right) ,\tau \right) \left( \det \frac{\partial 
\mathbf{Y}}{\partial \mathbf{y}}\right) ^{-1}\left( \det \frac{\partial X}{%
\partial x}\right) ^{-1}dv^{\prime }\mbox{d}X\mbox{d}\mathbf{Y}\leq 
\end{equation*}%
\begin{equation*}
\leq \left\vert V\right\vert e^{Ct}\int f^{p}\left( X\left( \tau \right)
,v^{\prime },\mathbf{Y}\left( \tau \right) ,\tau \right) \mbox{d}v^{\prime }%
\mbox{d}X\mbox{d}\mathbf{Y}.
\end{equation*}%
We take the $p$-th power of (\ref{fint}) and integrate over $x,$ $v,$ and $%
\mathbf{y}$ to get for $t\geq 0,$ 
\begin{equation}
\Vert f(t)\Vert \hbox{\raise -0.5mm \hbox{$_{L^{p}}$}}\leq \Vert f_{0}\Vert %
\hbox{\raise -0.5mm \hbox{$_{L^{p}}$}}+C\left( V\right) \!e^{Ct}\int_{0}^{t}%
\bigg\{(1+\Lambda \left( \left\Vert \mathbf{S}\left( \tau \right) \right\Vert
_{L^{\infty }}\right) +\Pi \left( \left\Vert \mathbf{S}\left( \tau \right)
\right\Vert _{L^{\infty }}\right) +  \label{fLp}
\end{equation}%
\begin{equation*}
+\left\Vert \frac{\partial \mathbf{S}}{\partial t}\left( \tau \right)
\right\Vert _{L^{\infty }}+\left\Vert \frac{\partial \mathbf{S}}{\partial x}%
\left( \tau \right) \right\Vert _{L^{\infty }}\bigg\}\times \Vert f\left(
\tau \right) \Vert \hbox{\raise -0.5mm
\hbox{$_{L^{p}}$}}\mbox{d}\tau .
\end{equation*}%
Using Lemma \ref{gradS}, we get for all $0\leq t\leq T$ and with $p=2$ 
\begin{equation}
\Vert f(t)\Vert \hbox{\raise -0.5mm \hbox{$_{L^2}$}}\leq \Vert f_{0}\Vert %
\hbox{\raise -0.5mm \hbox{$_{L^{2}}$}}+Ce^{Ct}\!\!\!\int_{0}^{t}\left[
1+\left\Vert n\left( 0\right) \right\Vert _{L^{1}}\left\{ 1+\ln \left(
\left\Vert n(\tau )\right\Vert _{L^{2}}+1\right) \right\} \right] \times
\Vert f(\tau )\Vert \hbox{\raise -0.5mm
\hbox{$_{L^2}$}}\}\mbox{d}\tau \leq   \label{f-L2!}
\end{equation}%
\begin{equation*}
\leq \Vert f_{0}\Vert \hbox{\raise -0.5mm \hbox{$_{L^{2}}$}}%
+Ce^{Ct}\!\!\!\int_{0}^{t}\left[ 1+\left\Vert f\left( 0\right) \right\Vert
_{L^{1}}\left\{ 1+\ln \left( \left\Vert f(\tau )\right\Vert
_{L^{2}}+1\right) \right\} \right] \times \Vert f(\tau )\Vert 
\hbox{\raise -0.5mm
\hbox{$_{L^2}$}}\}\mbox{d}\tau .
\end{equation*}%
By applying the Gronwall Lemma \ref{gronwall} to (\ref{f-L2!}), we obtain
for $t\geq 0,$ 
\begin{equation*}
\Vert f(t)\Vert \hbox{\raise -0.5mm \hbox{$_{L^2}$}}\leq C\left(
k_{i},k_{i}^{0},d_{i},\Lambda ,\Pi ,V,\Vert f_{0}\Vert 
\hbox{\raise -0.5mm
\hbox{$_{L^{1}}$}},\Vert f_{0}\Vert \hbox{\raise -0.5mm \hbox{$_{L^{2}}$}}%
\right) <\infty ,
\end{equation*}%
\begin{equation}
\left\Vert \frac{\partial \mathbf{S}}{\partial x}\left( t\right) \right\Vert
_{L^{\infty }}\leq C\left( k_{i},k_{i}^{0},d_{i},\Lambda ,\Pi ,V,\Vert
f_{0}\Vert \hbox{\raise -0.5mm \hbox{$_{L^{1}}$}},\Vert f_{0}\Vert %
\hbox{\raise -0.5mm \hbox{$_{L^{2}}$}}\right) <\infty,  \label{SLinf}
\end{equation}%
where we used Lemma \ref{gradS} to get estimate (\ref{SLinf}).
We now apply (\ref{SLinf}) to (\ref{fLp}) and we get (for $t\geq 0$ and for
all $1\leq p\leq \infty $)%
\begin{equation}
\Vert f(t)\Vert \hbox{\raise -0.5mm \hbox{$_{L^p}$}}\leq C\left(
k_{i},k_{i}^{0},d_{i},\Lambda ,\Pi ,V,\Vert f_{0}\Vert 
\hbox{\raise -0.5mm
\hbox{$_{L^{1}}$}},\Vert f_{0}\Vert \hbox{\raise -0.5mm \hbox{$_{L^{2}}$}}%
,\Vert f_{0}\Vert \hbox{\raise -0.5mm \hbox{$_{L^{2}}$}}\right) <\infty ,
\label{pestim}
\end{equation}%
i.e. we have obtained (\ref{fpestimate}). Using the elliptic equation (\ref%
{ellipticS}), the second derivative of the extracellular signal can be
expressed as 
\begin{equation}
\frac{\partial ^{2}S_{i}}{\partial x^{2}}=-\frac{k_{i}}{d_{i}}n+\frac{%
k_{i}^{0}}{d_{i}}S_{i}.  \label{boot}
\end{equation}%
Using (\ref{pestim}) and the elliptic theory, we deduce (\ref{Spestimate}).
Thus we complete the proof of Theorem \ref{theorem2}.

\rightline{Q.E.D.}

\medskip

\noindent We conclude this section with two corollaries. They provide
other conditions for the global existence of solutions. The proofs are
omitted because they are similar to proofs of Theorem \ref{theorem1} and
Theorem \ref{theorem2}.

\begin{corollary}
Assume $(\ref{growthylam})$ and $(\ref{growthT})$. Suppose there exists a
non-negative, nondecreasing continuous function $\Pi \left( \cdot \right)
\in C\left( \mathbb{R}\right) $ and $\gamma >0$ with $\omega \gamma \leq 1$
satisfying 
\begin{equation}
{\nabla }_{\mathbf{y}}{\cdot }\mathbf{F}\leq 0\text{ \ \ \ and \ \ }{\nabla }%
_{\mathbf{y}}{\cdot }\mathbf{F}\left( \mathbf{z},\mathbf{y}\right) \leq
C\left( 1+\Pi \left( \left\vert \mathbf{z}\right\vert \right) +\left\vert 
\mathbf{y}\right\vert ^{\gamma }\right)   \label{divF1}
\end{equation}%
Assume that $f_{0}\in L^{1}\cap L^{\infty }({}\mathbb{R}\times V\times 
\mathbb{R}^{m})$ and let the initial condition $S_{0}\in \lbrack W^{2,p}({}%
\mathbb{R})]^{M}$ satisfy $(\ref{ellipticS})$. Then there exists a global
solution of the system $(\ref{hyperpell})$ -- $(\ref{ellipticS})$
satisfying, for all $t\geq 0$ 
\begin{equation*}
f(\cdot ,\cdot ,\cdot ,t)\in L^{1}\cap L^{\infty }({}\mathbb{R}\times
V\times \mathbb{R}^{m}),
\end{equation*}%
\begin{equation*}
\mathbf{S}(\cdot ,t)\in \lbrack W^{2,p}(\mathbb{R})]^{M},\text{ \ \ \ \ for
all }1\leq p<+\infty 
\end{equation*}%
and initial conditions $f(\cdot ,\cdot ,\cdot ,0)=f_{0}(\cdot ,\cdot ,\cdot )
$ and $\mathbf{S}(\cdot ,0)=\mathbf{S}_{0}(\cdot )$. \label{cor1}
\end{corollary}

\begin{corollary}
Assume 
\begin{equation}
\lambda (y_{1})\leq C,\qquad T(v,v^{\prime },\mathbf{y})\leq C(1+|\lambda
(y_{1})|).  \label{growthTlam}
\end{equation}%
We further assume that $\mathbf{F}$ satisfies either $(\ref{divF2})$ or 
$(\ref{divF1})$. 
Assume that $f_{0}\in L^{1}\cap L^{\infty }({}\mathbb{R}\times
V\times \mathbb{R}^{m})$ and $\mathbf{S}_{0}\in \lbrack W^{1,\infty }({}%
\mathbb{R})]^{M}$ with compact support. Then there exists a global solution of
system of equations
$(\ref{hyperpell})$ and $(\ref{parabSsim})$ satisfying 
\begin{equation*}
f(\cdot ,\cdot ,\cdot ,t)\in L^{1}\cap L^{\infty }({}\mathbb{R}\times
V\times \mathbb{R}^{m}),
\end{equation*}%
\begin{equation*}
\mathbf{S}(\cdot ,t)\in \left[ W^{1,\infty }(\mathbb{R})\right] ^{M}
\end{equation*}%
and initial conditions $f(\cdot ,\cdot ,\cdot ,0)=f_{0}(\cdot ,\cdot ,\cdot )
$ and $\mathbf{S}(\cdot ,0)=\mathbf{S}_{0}(\cdot )$.
\end{corollary}

\section{Discussion}

\label{secdiscussion}

The simplified model of the bacterial signal transduction was studied 
in \cite{Erban:2004:ICB,Erban:2004:STS} where equation (\ref{rom14}) 
was given as (\ref{simmodel}). Using model (\ref{simmodel}) for the steady 
extracellular signal, one can derive the closed macroscopic (Keller-Segel,
chemotaxis) equation
for some parameter regimes. See \cite{Erban:2004:ICB} in 1D and \cite%
{Erban:2004:STS} in 2D/3D. Hence, the transport equation framework can be
used to study the macroscopic behaviour in terms of microscopic
parameters for the steady extracellular signals and simplified models 
of the signal transduction.

Here, we focused on more complex models where we coupled the complex
transport equation (\ref{hyperp}) with the parabolic or elliptic equation
for the signal (\ref{parabS}). The starting point of the analysis of such
complex models is the existence theory. In this paper, we provided several
sets of sufficient conditions for the global existence of solutions
of system (\ref{hyperp}) -- (\ref{parabS}). There 
are many open questions remaining, e.g. the existence theory
in $N$-dimensional physical space.  It is also not clear whether one can
derive the closed evolution equation for the density of cells $n(x,t)$ as
we did for the simple case of noninteracting particles \cite%
{Erban:2004:ICB,Erban:2004:STS}. If we are not able to derive 
the macroscopic equations then suitable computational approaches 
have to be used to study the macroscopic behaviour of bacteria
\cite{Erban:2005:CAE}.

There are several related results on kinetic models of the 
cellular movement. They often do not take the intracellular 
dynamics into account. Kinetic models were derived in \cite%
{Alt:1980:BRW,Othmer:1988:MDB} using stochastic models of
the movement of cells like bacteria or leukocytes. 
Reference \cite{Othmer:2002:DLT} addresses the formal
diffusion limit of kinetic models to the classical Keller-Segel model.
The discussion on issues of aggregation,
blow-up, and collapse for certain class of random walks can be 
found in \cite%
{Othmer:1997:ABC}. A Boltzmann-type kinetic model for chemotaxis 
 without the internal dynamics coupled with an elliptic equation for
the extracellular signal is studied in \cite{Chalub:2004:KMC} where 
global existence and rigorous diffusion limit to the Keller-Segel model 
were proven. In \cite{Hwang:2005:DDL,Hwang:2005:GSN}, a more general 
kinetic model was treated in
two and three dimensions. A one-dimensional hyperbolic model was
studied in 
\cite{Hwang:2005:GEC}. The papers 
\cite{Hwang:2005:DDL,Hwang:2005:GEC,Hwang:2005:GSN}
took into account the effect of the gradient
and the temporal derivative of the chemical signal and showed the global
existence of smooth solutions with smooth initial data as well as the
rigorous diffusive limit to the classical Keller-Segel model. However, all
the rigorous global existence results so far have not included the temporal
derivative of the signal in the growth condition of the turning frequency
as we did in this paper. See 
also \cite{Perthame:2004:PMC} for more related works.

\section*{Acknowledgements}

This work was partially supported by the Max Planck Institute for
Mathematics in Sciences, Biotechnology and Biological Sciences Research
Council, University of Oxford, Trinity College Dublin and Linacre College,
Oxford.

\providecommand{\bysame}{\leavevmode\hbox to3em{\hrulefill}\thinspace}


\end{document}